\documentclass[11pt]{article}
\usepackage[utf8]{inputenc}
\usepackage{amsbsy,amsfonts,amsmath,amssymb,amsthm
,bbm,enumerate,euscript,graphicx,color, subcaption}

\usepackage{hyperref} 
\hypersetup{
    colorlinks=true,
    linkcolor=blue,
    citecolor=blue,
    urlcolor = magenta,
}

\def\R{\mathbb{R}}
\def\N{\mathbb{N}}
\def\Z{\mathbb{Z}}
\def\E{\mathbb{E}}
\def\V{\mathrm{Var}}
\def\C{\mathrm{Cov}}

\def\I{\mathbf{1}}

\newtheorem{theo}{Theorem}

\newtheorem{prop}{Proposition}
\newtheorem{lemm}{Lemma}

\theoremstyle{definition}
\newtheorem{defn}{Definition}
\newtheorem*{proo}{Proof}
\newtheorem{exam}{Example}

\numberwithin{equation}{section}
\newtheorem*{rema}{Remark}

 \newtheorem{proposition}{Proposition} \newtheorem{remark}{Remark} 

\title{Decorated lattices: hyperuniformity and limit theorems} \author{Daniela Flimmel\thanks{Charles University, daniela.flimmel@matfyz.cuni.cz, corresponding author} 
\and Luca Lotz\thanks{German Aerospace Center and Karlsruhe Institute of Technology,
luca.lotz@dlr.de}}
\date{\today}

\begin{document}

\maketitle

\begin{abstract}
We introduce \emph{decorated lattices}, a point process obtained by placing finite random point configurations in lattice cells. A stationary mark field selects the motif laws and acts as a random environment transmitting dependence between cells, while allowing flexible local
arrangements. Using this construction, one can efficiently generate large regions of disordered patterns with desired local behavior while the random environment controls the asymptotic density fluctuations. 

Under moment and summability assumptions, we derive an explicit decomposition of the Bartlett spectral measure into contributions from conditional motif fluctuations and
correlations between marks.
This yields criteria for hyperuniformity, or generally $p$-uniformity. Independent, nondegenerate count fluctuations destroy hyperuniformity, whereas suitable dependence governed by the random environment can restore it despite random and unequal cell counts.
Once count fluctuations are sufficiently suppressed, motif geometry and its spatial correlations govern the finer order of uniformity.
Under mixing assumptions, we also establish central limit theorems for linear statistics, with fluctuations governed by cell counts or, at a finer scale, by the first spatial
moments of the motifs.
\end{abstract}

\vspace{0.2cm}



\section{Introduction}

Point processes provide a natural framework for modelling spatial configurations of particles or the locations of constituents of a heterogeneous material; see
\cite{LP17,DVJ08} for general background. In such applications, both the local geometry and the fluctuations of the number of constituents over large regions can be relevant. Lattice-based constructions offer a convenient way to prescribe local configurations while retaining a tractable description of their spatial dependence.
A central question is how much freedom these constructions allow at the local scale without losing control of their large-scale fluctuations.

This question is closely connected to hyperuniformity. A stationary point process $\Xi$ is called hyperuniform if
\[
\frac{\V(\Xi(B_r))}{\lambda_d(B_r)}
\longrightarrow 0
\qquad\text{as }r\to\infty,
\]
where $B_r$ is the ball of radius $r$ centred at the origin. Thus, its number variance grows more slowly than the volume of the observation window, in contrast to the volume-order fluctuations of a homogeneous Poisson point process.
The concept was developed in statistical physics in \cite{TS03} and provides a common description of suppressed density fluctuations in crystals, many quasicrystals, and certain disordered systems. For perspectives from physics and materials science, see \cite{T18}; mathematical surveys include \cite{GL17,C21, LR25}.

The coexistence of local disorder and suppressed density fluctuations has motivated applications to materials design. For example, \cite{FTS09} constructs disordered dielectric
structures with large, complete photonic band gaps from hyperuniform point patterns.
A related construction principle appears in  \cite{KT19}, where space is partitioned
into cells and particles are placed within them so that each cell has the same local packing fraction.
This produces hyperuniform two-phase materials and illustrates how local constraints can control fluctuations at much larger scales. 

Beyond identifying hyperuniform models, an important practical question is how to generate large realizations efficiently. This challenge is highlighed in 
\cite[Section~3.6]{C21}, with particular emphasis on amorphous
hyperuniform patterns containing millions of points.
Perturbed lattices provide a computationally simple starting
point (\cite{DFHL24,F26}), but their spectra generally retain Bragg peaks inherited
from the underlying lattice. These peaks can nevertheless be
eliminated by a suitable choice of the perturbation law:
the cloaking construction of 
\cite{KKT20} produces large hyperuniform samples without Bragg
peaks, although the lattice structure may remain detectable
through higher-order correlations.
More recently, \cite{KLLY25} establishes conditions under
which invariant transports preserve hyperuniformity and
constructs hyperuniform point processes by sampling within
the cells of fair partitions.
The work \cite{LK26} develops these ideas further through
spatial moment constraints, including constructions of
isotropic processes with arbitrarily high orders of
$p$-uniformity that admit simulation in linear time.

We investigate a complementary mechanism: controlling
density fluctuations through dependence between cells,
while allowing their populations and internal arrangements
to vary. We introduce a stationary random environment that selects
the motif laws and transmits dependence between cells.
This makes it possible to prescribe local arrangements
while controlling large-scale count fluctuations through
the environment. The resulting model also allows for fast generation of large datasets with desired level of $p$-uniformity.

In this paper, we introduce \emph{decorated lattices} driven by a
stationary random environment. Each cell $k+Q$, where $Q=[0,1)^d$, contains a finite random
point configuration $\Xi_k$, called a motif.
Conditionally on $\mathbf M$, the motifs are independent,
with $\Xi_k$ having law $P_{M_k}$.
An independent uniform shift $U$ in $Q$ gives the stationary
point process
\[
\widetilde\Phi
=
\sum_{k\in\mathbb Z^d}\theta_{k+U}\Xi_k.
\]
The marks may specify different local patterns, sampling
distributions, or numbers of points. Dependence between cells
is transmitted through the mark field, while the points
within an individual motif may have arbitrary dependence.
The construction includes for instance perturbed lattices and homogeneous
Poisson point processes, as well as their combinations. It also separates simulation into sampling the environment
and then generating the motifs independently given that environment.

Our construction is related to the independent cluster
perturbations studied in \cite[Proposition~5]{La25}, where
independent copies of a finite random configuration are
attached to the points of a stationary parent process.
Here, we confine the motifs to bounded lattice cells,
but allow dependence between them through the random
environment. In comparison with the fixed-count
partition constructions discussed above, we allow
random and unequal cell counts whose correlations
can suppress their collective fluctuations.

Our first main result is an explicit decomposition of the
Bartlett spectral measure of $\widetilde\Phi$ under a
second-moment condition and a summability assumption on the dependence between marks. The spectrum consists of an absolutely continuous part and Bragg peaks at nonzero reciprocal-lattice points. Its continuous density separates the fluctuations remaining within a cell after its mark has been specified from the correlations transmitted by
the environment. This decomposition
makes it possible to identify how the different ingredients
of the construction contribute to the spectrum.

In particular, the behaviour at frequency zero is governed
by the cell counts
\[
N_k:=\Xi_k(Q), \quad k \in \Z^d.
\]
Writing $S$ for the density of the absolutely continuous
part of the Bartlett spectral measure, we obtain
\[
S(0)
=
\E[\V(N_0\mid M_0)]
+
\sum_{k\in\mathbb Z^d}
\C\left(
\E[N_0\mid M_0],
\E[N_k\mid M_k]
\right).
\]
Both terms are nonnegative. Consequently, under our assumptions, $\widetilde\Phi$ is hyperuniform
if and only if the counts are determined by the marks
and the field $(\E[N_k\mid M_k])_{k\in\mathbb Z^d}$
is hyperuniform. Equivalently, hyperuniformity of the
decorated lattice is characterized by hyperuniformity
of its cell-count field. 

\begin{figure}[t]
    \centering
    \includegraphics[width=\textwidth]{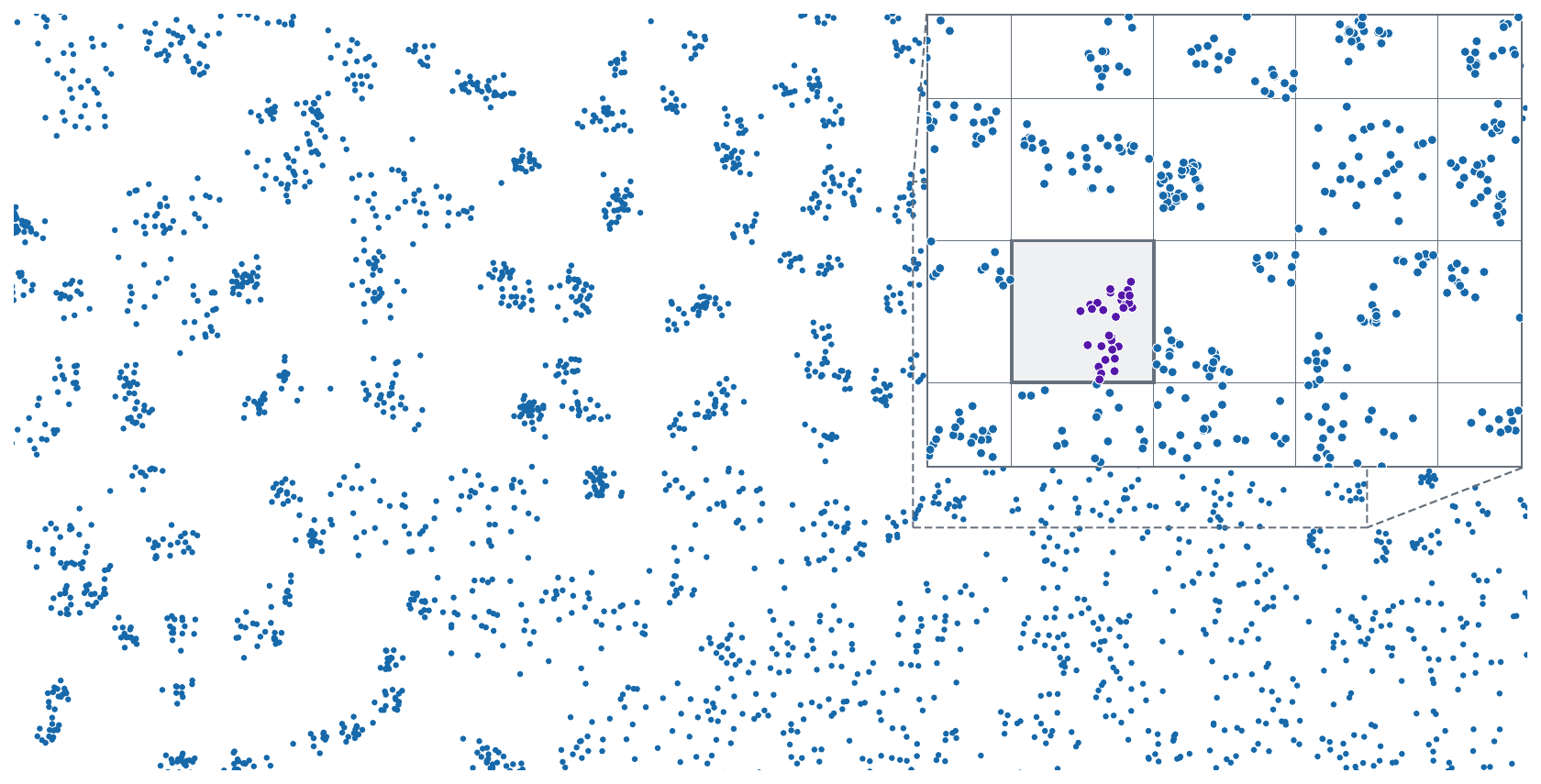}
    \caption{Sample of a stationarized decorated lattice driven by a stationary random environment. Its mark field has two spatially correlated components. The first component determines how much the points inside each cell cluster together, and the second governs the number of clusters per cell. All additional randomness is conditionally independent given the mark field. In this example, we placed exactly $30$ points in each cell. Hence, Proposition~\ref{P_hyperuniformity_decorated} yields that the stationarized decorated lattice is hyperuniform. In particular, it is $2$-uniform and not beyond \linebreak$2$-uniform.}
\end{figure}

This criterion provides a construction principle in which the random environment controls hyperuniformity.
There is no requirement that every cell contain the same
number of points: random, unequal counts are compatible
with hyperuniformity when their collective fluctuations
are sufficiently suppressed. In contrast, residual
conditional randomness in the counts contributes a positive
amount to $S(0)$ that cannot be cancelled by correlations
between marks. Subject to the stated assumptions, the
positions of points within the motifs do not affect this
criterion, leaving considerable freedom to prescribe
local arrangements.

Local geometry becomes relevant when one asks how rapidly
fluctuations vanish. We investigate this finer behaviour
through the notion of $p$-uniformity introduced in \cite{LK26},
which quantifies the spectral mass near the origin.
Our spectral formulas yield sufficient conditions for
$p$-uniformity and allow explicit calculations in examples.
These include binomial motifs, correlated cell counts,
random within-cell sampling distributions, and motifs
whose spatial moments satisfy deterministic constraints.
The examples distinguish the role of count correlations
in establishing hyperuniformity from that of spatial
moment constraints in producing higher orders of uniformity.

Our second main contribution concerns the corresponding
fluctuation limits.
Under polynomial $\alpha$-mixing and moment assumptions,
we establish central limit theorems for centred linear
statistics at two scales.
At normalization $r^{-d/2}$, the limiting variance is
determined by the summed covariances of the cell-count
field.
When that field is beyond $2$-uniform, its centred
contribution becomes negligible at the finer normalization
$r^{1-d/2}$.
The fluctuations are then governed by the first spatial
moments
\[
T_k:=\int_Q x\,\Xi_k(\mathrm dx),
\qquad k\in\mathbb Z^d,
\]
with limiting variance expressed through the gradient
of the test function and the summed covariances of
$(T_k)_{k\in\mathbb Z^d}$.
The proofs combine the transport expansion of \cite{LK26}
with the weighted random-field central limit theorem
of \cite{JP09}.
Related results for independently perturbed lattices
include the variance asymptotics for Gaussian perturbations
in \cite{Ya21} and the fluctuation limits for more general
independent perturbations in \cite{Ma25}.
Our results address bounded motifs with random cell
populations and dependence transmitted through a stationary
environment, identifying how this dependence enters both
the hyperuniformity criterion and the limiting variances.

The paper is organized as follows.
Section~\ref{sec:2} collects the preliminaries on Fourier
transforms, point processes, and hyperuniformity.
Section~\ref{sec:3} introduces decorated lattices driven
by a stationary random environment.
Section~\ref{sec:4} derives their intensity and second-order
measures, establishes the spectral decomposition, and
develops criteria for hyperuniformity and $p$-uniformity.
Section~\ref{sec:5} proves central limit theorems at the
scales governed by cell counts and first spatial moments.
Finally, Section~\ref{sec:6} presents explicit constructions
illustrating the effects of count dependence and spatial moment constraints. 

\section{Preliminaries}\label{sec:2}

\subsection{Fourier transform of positive definite measures}
For background on Fourier analysis and distributions, we refer to
\cite[Chapters~8--9]{Folland1999}. We specify our Fourier normalization
below.

For $f\in L^1(\mathbb R^d)$, we use the Fourier transform
\[
\widehat f(k)
:=
\int_{\mathbb R^d}
e^{-ik\cdot x}f(x)\,\mathrm dx,
\qquad k\in\mathbb R^d.
\]
For $h\in L^1(\mathbb R^d)$, define the inverse Fourier transform by
\[
\check h(x)
:=
\frac{1}{(2\pi)^d}
\int_{\mathbb R^d}
h(k)e^{ik\cdot x}\,\mathrm dk.
\]
In particular, Fourier inversion gives $\check{\hat f}=f$
for $f\in\mathcal S(\mathbb R^d)$, where $\mathcal S(\mathbb R^d)$
denotes the Schwartz space.

For $f,g\in\mathcal S(\mathbb R^d)$, define the tilted convolution by
\[
(f\star g)(x)
:=
\int_{\mathbb R^d}
f(y)\overline{g(y-x)}\,\mathrm dy.
\]
Equivalently,
\[
f\star g=f*\widetilde g,
\qquad
\widetilde g(x):=\overline{g(-x)},
\]
where $*$ denotes ordinary convolution. With the Fourier convention
above,
\[
\widehat{f\star g}(k)
=
\widehat f(k)\overline{\widehat g(k)}
\]
and
\[
\widehat{fg}
=
\frac{1}{(2\pi)^d}\,
\widehat f*\widehat g.
\]

A locally finite real signed Radon measure $\nu$ on $\mathbb R^d$
is called \emph{positive definite} if
\[
\nu(f\star f)\geq 0
\]
for every $f\in C_c^\infty(\mathbb R^d;\mathbb C)$, the space of
smooth, compactly supported complex-valued functions.

By the Bochner--Schwartz theorem, a positive-definite measure
defines a tempered distribution and has a unique positive
tempered spectral measure $\widehat\nu$. With our convention,
this measure is characterized by
\[
\nu(f\star g)
=
\frac{1}{(2\pi)^d}
\int_{\mathbb R^d}
\widehat f(k)\overline{\widehat g(k)}
\,\widehat\nu(\mathrm dk),
\qquad
f,g\in C_c^\infty(\mathbb R^d;\mathbb C).
\]
In particular, $\widehat\nu$ is locally finite.
If, additionally, $\nu$ has finite total variation, define its
Fourier--Stieltjes transform by
\[
\mathcal F(\nu)(k)
:=
\int_{\mathbb R^d}
e^{-ik\cdot x}\,\nu(\mathrm dx).
\]
Fourier inversion and Fubini's theorem give
\[
\nu(f\star g)
=
\frac{1}{(2\pi)^d}
\int_{\mathbb R^d}
\widehat f(k)\overline{\widehat g(k)}
\,\mathcal F(\nu)(-k)\,\mathrm dk.
\]
Consequently, $\widehat\nu$ is absolutely continuous with respect
to Lebesgue measure. Since $\nu$ is real and positive definite, it is symmetric
under reflection, 
\[
\widehat\nu(\mathrm dk)
=
\mathcal F(\nu)(-k)\,\mathrm dk
=
\mathcal F(\nu)(k)\,\mathrm dk.
\]
Its density $\mathcal F(\nu)$ is continuous, even, real-valued,
and nonnegative.

\subsection{Point processes and second-order measures}

We regard point processes as random locally finite counting
measures on $\mathbb R^d$; see \cite{LP17,DVJ08}.
More precisely, a point process $\Xi$ is a random element of
$\mathcal N(\mathbb R^d)$, equipped with the $\sigma$-algebra
generated by the evaluation maps $\nu\mapsto\nu(B)$ for bounded
Borel sets $B$.
The count $\Xi(B)$ includes multiplicities.
The process is \emph{simple} if, almost surely,
$\Xi(\{x\})\leq1$ for every $x\in\mathbb R^d$;
simplicity is not assumed unless stated explicitly.

For bounded, compactly supported Borel functions
$f:\mathbb R^d\to\mathbb C$, write
\[
\Xi(f):=\int_{\mathbb R^d}f(x)\,\Xi(\mathrm dx).
\]
The \emph{intensity measure} is
$\Lambda_\Xi(B):=\E[\Xi(B)]$.
With $\theta_u\nu(B):=\nu(B-u)$, the process is
\emph{stationary} if $\theta_u\Xi\overset{d}{=}\Xi$ for every
$u\in\mathbb R^d$.
If $\Lambda_\Xi$ is locally finite, stationarity implies
\[
\Lambda_\Xi=\gamma\lambda_d,
\qquad
\gamma:=\E[\Xi([0,1)^d)],
\]
where $\lambda_d$ is Lebesgue measure and $\gamma$ is the
\emph{intensity}.

\paragraph{Moment measures.}

For the second-order considerations below, assume that $\Xi$
is \emph{locally square-integrable}:
\[
\E[\Xi(B)^2]<\infty
\qquad\text{for every bounded Borel set }B.
\]
The \emph{second-order moment measure}
$\alpha_\Xi^{(2)}:=\E[\Xi\otimes\Xi]$ is locally finite and
uniquely determined by
\[
\alpha_\Xi^{(2)}(C\times D)=\E[\Xi(C)\Xi(D)]
\]
for bounded Borel sets $C,D$.
For every nonnegative measurable function $F$,
\[
\alpha_\Xi^{(2)}(F)
=
\E\left[
\iint F(x,y)\,\Xi(\mathrm dx)\Xi(\mathrm dy)
\right].
\]

Let $\Delta_{\mathbb R^d}(x):=(x,x)$.
Its pushforward acts on measures by
\[
\bigl((\Delta_{\mathbb R^d})_\#\mu\bigr)(C\times D)
=
\mu(C\cap D).
\]
The \emph{second-order factorial moment measure}
\[
\alpha_{\Xi,!}^{(2)}
:=
\alpha_\Xi^{(2)}
-
(\Delta_{\mathbb R^d})_\#\Lambda_\Xi
\]
therefore satisfies
\[
\alpha_{\Xi,!}^{(2)}(C\times D)
=
\E\bigl[\Xi(C)\Xi(D)-\Xi(C\cap D)\bigr].
\]
It counts ordered pairs of distinct point occurrences,
which may share a location if $\Xi$ is not simple.

\paragraph{Covariance and cumulant measures.}

The \emph{second-order cumulant measure}, or
\emph{covariance measure}, and the
\emph{second-order factorial cumulant measure} are
\[
\beta_\Xi^{(2)}
:=
\alpha_\Xi^{(2)}-\Lambda_\Xi\otimes\Lambda_\Xi,
\qquad
\beta_{\Xi,!}^{(2)}
:=
\alpha_{\Xi,!}^{(2)}-\Lambda_\Xi\otimes\Lambda_\Xi.
\]
They are locally finite signed measures satisfying
\[
\beta_\Xi^{(2)}(C\times D)
=
\C\bigl(\Xi(C),\Xi(D)\bigr),
\qquad
\beta_\Xi^{(2)}
=
\beta_{\Xi,!}^{(2)}
+
(\Delta_{\mathbb R^d})_\#\Lambda_\Xi.
\]
For complex-valued random variables $X, Y$, we use
\[
\C(X,Y):=\E[(X-\E X)\overline{(Y-\E Y)}],
\qquad
\V(X):=\C(X,X).
\]

\paragraph{Reduced second-order measures.}

Suppose additionally that $\Xi$ is stationary, with intensity
$\gamma$, and put $Q:=[0,1)^d$.
Stationarity allows the second-order measures to be expressed
in terms of relative displacements.
We use the convention $x-y$ and do not normalize the reduced
measures by the intensity; see \cite{DVJ08}.

The \emph{reduced second-order moment measure} is
\begin{equation}\label{EQ_reduced_second_moment}
\alpha_{\Xi,\mathrm{red}}^{(2)}(B)
:=
\int_{Q\times\mathbb R^d}
\mathbf 1_B(x-y)\,
\alpha_\Xi^{(2)}(\mathrm dx,\mathrm dy),
\qquad B\in\mathcal B(\mathbb R^d).
\end{equation}
Similarly, its factorial counterpart is
\[
\alpha_{\Xi,\mathrm{red},!}^{(2)}(B)
:=
\int_{Q\times\mathbb R^d}
\mathbf 1_B(x-y)\,
\alpha_{\Xi,!}^{(2)}(\mathrm dx,\mathrm dy).
\]
Both measures are locally finite.
By stationarity, replacing $Q$ by any bounded Borel set
of unit volume leaves these definitions unchanged.

The \emph{reduced covariance measure} and
\emph{reduced factorial cumulant measure} are
\[
\beta_{\Xi,\mathrm{red}}^{(2)}
:=
\alpha_{\Xi,\mathrm{red}}^{(2)}-\gamma^2\lambda_d,
\qquad
\beta_{\Xi,\mathrm{red},!}^{(2)}
:=
\alpha_{\Xi,\mathrm{red},!}^{(2)}-\gamma^2\lambda_d.
\]
Since self-pairs have displacement zero,
\[
\alpha_{\Xi,\mathrm{red}}^{(2)}
=
\alpha_{\Xi,\mathrm{red},!}^{(2)}+\gamma\delta_0,
\qquad
\beta_{\Xi,\mathrm{red}}^{(2)}
=
\beta_{\Xi,\mathrm{red},!}^{(2)}+\gamma\delta_0.
\]

For bounded, compactly supported Borel functions $f,g$,
stationarity gives
\[
\alpha_{\Xi,\mathrm{red}}^{(2)}(f\star g)
=
\E\bigl[\Xi(f)\overline{\Xi(g)}\bigr],
\qquad
\beta_{\Xi,\mathrm{red}}^{(2)}(f\star g)
=
\C\bigl(\Xi(f),\Xi(g)\bigr).
\]
Taking $g=f$ shows that both
$\alpha_{\Xi,\mathrm{red}}^{(2)}$ and
$\beta_{\Xi,\mathrm{red}}^{(2)}$ are positive definite.

\paragraph{Bartlett spectral measure.}

The \emph{Bartlett spectral measure} of $\Xi$ is
\[
\widehat\beta_\Xi
:=
\widehat{\beta_{\Xi,\mathrm{red}}^{(2)}}.
\]
It is a positive, locally finite measure satisfying
\[
\C\bigl(\Xi(f),\Xi(g)\bigr)
=
\frac{1}{(2\pi)^d}
\int_{\mathbb R^d}
\widehat f(k)\overline{\widehat g(k)}
\,\widehat\beta_\Xi(\mathrm dk),
\qquad
f,g\in C_c^\infty(\mathbb R^d;\mathbb C).
\]

Likewise, write
$\widehat\alpha_\Xi:=
\widehat{\alpha_{\Xi,\mathrm{red}}^{(2)}}$.
Then
\[
\E\bigl[\Xi(f)\overline{\Xi(g)}\bigr]
=
\frac{1}{(2\pi)^d}
\int_{\mathbb R^d}
\widehat f(k)\overline{\widehat g(k)}
\,\widehat\alpha_\Xi(\mathrm dk).
\]
The two spectral measures are related by
\[
\widehat\alpha_\Xi
=
\widehat\beta_\Xi+(2\pi)^d\gamma^2\delta_0.
\]
Thus, centering removes the contribution of the mean at
frequency zero. The covariance spectrum may still have
an atom at zero.

If
$\widehat\beta_\Xi(\mathrm dk)=S_\Xi(k)\,\mathrm dk$,
then $S_\Xi$ is called the \emph{spectral density}.
For positive intensity $\gamma$, the normalized density
$\gamma^{-1}S_\Xi$ is the \emph{structure factor}.

\begin{exam}[Poisson point process]
For a homogeneous Poisson point process $\Pi$ of intensity
$\gamma>0$,
\[
\alpha_{\Pi,!}^{(2)}
=
\gamma^2\lambda_d\otimes\lambda_d,
\qquad
\beta_{\Pi,!}^{(2)}=0,
\qquad
\beta_{\Pi,\mathrm{red}}^{(2)}=\gamma\delta_0.
\]
Hence,
\begin{equation}\label{Eq_Bartlett_Poisson}
\widehat\beta_\Pi=\gamma\lambda_d.
\end{equation}
Its spectral density is $S_\Pi(k)=\gamma$, and its structure
factor is identically one.
\end{exam}

\subsection{Hyperuniformity}

Hyperuniformity describes the suppression of large-scale number
fluctuations relative to those of a homogeneous Poisson process;
see \cite{TS03,T18}.
Throughout this subsection, let $B_r$ denote the Euclidean ball
of radius $r>0$ centred at the origin.

\begin{defn}[Hyperuniform point process]
A stationary, locally square-integrable point process $\Xi$
on $\mathbb R^d$ is called \emph{hyperuniform} if
\[
\frac{\V(\Xi(B_r))}{\lambda_d(B_r)}
\xrightarrow{r\to\infty}0.
\]
\end{defn}

Thus, the number variance grows sublinearly in the volume of
the observation window. In comparison, a homogeneous Poisson
process $\Pi$ with intensity $\gamma>0$ satisfies
\[
\V(\Pi(B_r))=\gamma\lambda_d(B_r).
\]

The connection with the spectral description follows from
\[
\V(\Xi(B_r))
=
\frac{1}{(2\pi)^d}
\int_{\mathbb R^d}
\left|\widehat{\mathbf 1_{B_r}}(k)\right|^2
\,\widehat\beta_\Xi(\mathrm dk).
\]
For Euclidean balls, hyperuniformity is equivalent to
\[
\frac{\widehat\beta_\Xi(B_\varepsilon)}{\varepsilon^d}
\xrightarrow{\varepsilon\to 0}0;
\]
see \cite{BH24}.
This formulation does not require the Bartlett spectral measure
to be absolutely continuous and therefore also applies to
processes with Bragg peaks or a singular continuous spectrum.

If $\widehat\beta_\Xi$ is absolutely continuous in a neighbourhood
of the origin and its density $S_\Xi$ is continuous at the origin,
then
\[
\lim_{\varepsilon\to 0}
\frac{\widehat\beta_\Xi(B_\varepsilon)}{\varepsilon^d}
=
\lambda_d(B_1)S_\Xi(0).
\]
Under these assumptions, hyperuniformity is therefore equivalent
to $S_\Xi(0)=0$.

To quantify the rate at which the spectral mass vanishes near
the origin, we use the notions of $p$-uniformity and beyond
$p$-uniformity introduced in \cite{LK26}.
We restrict the definition here only to point processes.

\begin{defn}[$p$-uniformity]\label{d:p-uniformity}
Let $p\in[-d,\infty)$, and let $\Xi$ be a stationary, locally
square-integrable point process on $\mathbb R^d$.
The process $\Xi$ is called \emph{$p$-uniform} if
\[
\frac{\widehat\beta_\Xi(B_\varepsilon)}{\varepsilon^d}
=
O(\varepsilon^p),
\qquad \varepsilon\to 0.
\]
It is called \emph{beyond $p$-uniform} if
\[
\frac{\widehat\beta_\Xi(B_\varepsilon)}{\varepsilon^d}
=
o(\varepsilon^p),
\qquad \varepsilon\to 0.
\]
\end{defn}

In particular, hyperuniformity is equivalent to beyond
$0$-uniformity, and $p$-uniformity for any $p>0$ implies
hyperuniformity. A homogeneous Poisson process of positive
intensity is $0$-uniform but not beyond $0$-uniform.

\section{Decorated lattices in random environments}\label{sec:3}

We construct a point process by placing a finite point configuration,
called a \emph{motif}, in each cell of the integer lattice.
The distribution of the motif in a given cell is selected by a
stationary random mark field. Conditional on this field, the motifs
are sampled independently. Thus, the marks describe a random
environment that introduces dependence between cells, while the
motif laws determine the local point configurations.

Let $Q:=[0,1)^d$ be the reference cell, and let
$\mathcal N_f(Q)$ denote the space of finite counting measures on $Q$,
equipped with the $\sigma$-algebra generated by the evaluation maps
$\nu\mapsto\nu(B)$, $B\in\mathcal B(Q)$.

\begin{defn}[Decorated lattice]
Let $\mathbf M=(M_k)_{k\in\mathbb Z^d}$ be a stationary random field
taking values in a standard Borel space $(\mathbb M,\mathcal M)$.
For each mark $m\in\mathbb M$, let $P_m$ be a probability law on
$\mathcal N_f(Q)$, such that $m\mapsto P_m(A)$ is measurable for every measurable $A \subset \mathcal{N}_f(Q)$.

Conditionally on $\mathbf M$, let the motifs
$(\Xi_k)_{k\in\mathbb Z^d}$ be independent with
\[
\mathcal L(\Xi_k\mid\mathbf M)=P_{M_k},
\qquad k\in\mathbb Z^d.
\]
The associated \emph{decorated lattice} is
\[
\Phi
:=
\sum_{k\in\mathbb Z^d}\theta_k\Xi_k,
\]
where
\[
\theta_u\nu(B):=\nu(B-u),
\qquad
u\in\mathbb R^d,\quad B\in\mathcal B(\mathbb R^d).
\]
We consider here its stationary version defined by
\[
\widetilde\Phi
:=
\theta_U\Phi
=
\sum_{k\in\mathbb Z^d}\theta_{k+U}\Xi_k,
\]
where $U$ is uniformly distributed on $Q$ and independent of
$\mathbf M$ and $(\Xi_k)_{k\in\mathbb Z^d}$.
\end{defn}

\begin{rema}
The conditional sampling rule can equivalently be written as
\[
\mathcal L\bigl((\Xi_k)_{k\in\mathbb Z^d}\mid\mathbf M\bigr)
=
\bigotimes_{k\in\mathbb Z^d}P_{M_k}
\qquad\text{almost surely}.
\]
\end{rema}

Each motif $\Xi_k$ is described in the reference cell $Q$ and then
translated into $k+Q$. Since a bounded set intersects only finitely
many lattice cells and every motif is finite, $\Phi$ and
$\widetilde\Phi$ are almost surely locally finite.
No moment assumptions on the motif sizes are needed for this
construction.

Stationarity of the mark field and the common sampling rule imply
that the motif field is stationary under lattice translations.
Consequently, $\Phi$ is invariant in distribution under translations
by vectors in $\mathbb Z^d$. The independent uniform shift $U$
makes $\widetilde\Phi$ stationary under all translations of
$\mathbb R^d$. 

The construction includes both deterministic and random motifs.
For example, if
\[
P_m=\delta_{\nu_m}
\]
for a measurable family of finite counting measures
$(\nu_m)_{m\in\mathbb M}$, then $\Xi_k=\nu_{M_k}$ almost surely. In this case,
the mark completely determines the motif.
More generally, $P_m$ may allow randomness in both the number
and the positions of the points.
Although motifs in distinct cells are conditionally independent,
they need not be independent unconditionally, because their marks
may be dependent.

For later use, write
\[
\pi:=\mathcal L(M_0),
\qquad
\pi_r:=\mathcal L(M_0,M_r),
\qquad r\in\mathbb Z^d.
\]
By stationarity, both marginals of $\pi_r$ equal $\pi$.
The joint law of the mark and motif in a typical cell is
\[
\varrho(\mathrm dm,\mathrm d\nu)
=
\pi(\mathrm dm)P_m(\mathrm d\nu),
\]
and the unconditional motif law is the mixture
\[
\overline P(A)
=
\int_{\mathbb M}P_m(A)\,\pi(\mathrm dm),
\qquad A\in\mathcal B(\mathcal N_f(Q)).
\]

The construction accommodates a broad range of spatial patterns,
from perturbed lattices with bounded perturbations to the
homogeneous Poisson point process. By varying the mark field and
the motif laws, one can control both the local arrangement of points
and the dependence between different regions. 
\begin{exam}[Bounded perturbation of a lattice]
Suppose that, for $\pi$-almost every $m$, the law $P_m$ is
concentrated on motifs containing exactly one point.
Then
\[
\Xi_k=\delta_{Y_k},
\qquad Y_k\in Q,
\]
where the random vectors $(Y_k)_{k\in\mathbb Z^d}$ are
conditionally independent given $\mathbf M$.
Thus,
\[
\Phi=\sum_{k\in\mathbb Z^d}\delta_{k+Y_k},
\qquad
\widetilde\Phi
=
\sum_{k\in\mathbb Z^d}\delta_{k+U+Y_k}.
\]
This is a perturbed lattice with a stationary, uniformly bounded
displacement field $(Y_k)_{k\in\mathbb Z^d}$.
The displacements may be dependent through the mark field.
In this case, $\widetilde{\Phi}$ is a simple point processes of unit intensity.
\end{exam}

\begin{exam}[Homogeneous Poisson point process]
Let $\gamma>0$, take a singleton mark space $\mathbb M=\{m\}$,
and let $P_m$ be the law of a Poisson point process on $Q$
with intensity measure $\gamma\lambda_d|_Q$.
Equivalently, each motif consists of a Poisson-distributed number
of points with mean $\gamma$, placed independently and uniformly
in $Q$, independently across cells.

Since the cells $(k+Q)_{k\in\mathbb Z^d}$ partition $\mathbb R^d$,
the superposition
\[
\Phi=\sum_{k\in\mathbb Z^d}\theta_k\Xi_k
\]
is a homogeneous Poisson point process with intensity $\gamma$.
Indeed, counts in disjoint bounded Borel sets are independent
Poisson random variables with means given by $\gamma\lambda_d$.
Thus, $\Phi$ is already stationary under all translations,
and the independent uniform shift does not change its law:
$\widetilde\Phi\overset{d}{=}\Phi$.
\end{exam}

\section{Bartlett spectrum and hyperuniformity}\label{sec:4}

In this section, we study the stationary decorated lattice $\widetilde\Phi$.

\subsection{Intensity measure}

For each $m\in\mathbb M$, define the motif intensity measure by
\[
\mu_m(C)
:=
\int_{\mathcal N_f(Q)}
\nu(C)\,P_m(\mathrm d\nu),
\qquad C\in\mathcal B(Q),
\]
and write
\[
n(m):=\mu_m(Q)
\]
for the expected number of points in a motif of type $m$.
Note that these quantities may in general be infinite. For every $C\in\mathcal B(Q)$, the map $m\mapsto\mu_m(C)$
is measurable, since $\nu\mapsto\nu(C)$ is measurable and
$m\mapsto P_m$ is a probability kernel.
In particular, $m\mapsto n(m)$ is measurable. Define the averaged motif intensity measure by
\[
\overline\mu(C)
:=
\int_{\mathbb M}\mu_m(C)\,\pi(\mathrm dm),
\qquad C\in\mathcal B(Q),
\]
and let
\[
\overline n
:=
\overline\mu(Q)
=
\int_{\mathbb M}n(m)\,\pi(\mathrm dm).
\]
In terms of the unconditional motif law $\overline P$, this gives
\[
\overline\mu(C)
=
\int_{\mathcal N_f(Q)}
\nu(C)\,\overline P(\mathrm d\nu).
\]
Thus, $\overline n$ is the expected number of points in a typical
cell. If $\overline n<\infty$, then $\overline\mu$ is finite,
and $\mu_m$ is finite for $\pi$-almost every $m$. The independent uniform shift spreads the averaged motif
intensity over space. Since each lattice cell has unit volume,
the resulting intensity equals $\overline n$.

\begin{prop}[Intensity of the stationary decorated lattice]
Assume that $\overline n<\infty$.
Then $\widetilde\Phi$ has intensity measure
\[
\Lambda_{\widetilde\Phi}
=
\overline n\,\lambda_d.
\]
\end{prop}

\begin{proo}
By the conditional sampling rule, for every $k\in\mathbb Z^d$
and $C\in\mathcal B(Q)$,
\[
\E[\Xi_k(C)\mid\mathbf M]
=
\int_{\mathcal N_f(Q)}
\nu(C)\,P_{M_k}(\mathrm d\nu)
=
\mu_{M_k}(C)
\qquad\text{a.s}.
\]
Taking expectations and using stationarity of the mark field,
we obtain
\begin{equation}\label{EQ_unconditional_first_moment_cell}
\E[\Xi_k(C)]
=
\int_{\mathbb M}\mu_m(C)\,\pi(\mathrm dm)
=
\overline\mu(C).
\end{equation}
Let $B\subseteq\mathbb R^d$ be a bounded Borel set.
By Tonelli's theorem, independence of $U$ from the motif field,
and \eqref{EQ_unconditional_first_moment_cell},
\begin{align*}
\Lambda_{\widetilde\Phi}(B)
&=
\E\left[
\sum_{k\in\mathbb Z^d}
\int_Q
\mathbf 1_B(k+x+U)\,\Xi_k(\mathrm dx)
\right]
=
\sum_{k\in\mathbb Z^d}
\int_Q\int_Q
\mathbf 1_B(k+x+u)\,
\mathrm du\,\overline\mu(\mathrm dx)
\\
&=
\int_Q
\left(
\sum_{k\in\mathbb Z^d}
\int_{k+x+Q}\mathbf 1_B(y)\,\mathrm dy
\right)
\overline\mu(\mathrm dx)
=
\lambda_d(B)\,\overline\mu(Q)
=
\overline n\,\lambda_d(B).
\end{align*}
Here we used that, for each fixed $x\in Q$, the sets
$(k+x+Q)_{k\in\mathbb Z^d}$ form a measurable partition
of $\mathbb R^d$.
The equality extends to all Borel sets by monotone convergence.
\qed
\end{proo}

\begin{rema}
Conditional independence of the motifs is not needed for this
first-order result. The proof uses only the common mean measure
$\E[\Xi_k(\cdot)]=\overline\mu$ and the independent uniform shift.
In particular, neither correlations between marks nor the
arrangement of points within motifs affect the intensity,
provided that $\overline n$ remains unchanged.
\end{rema}

\subsection{Second-order measures}

We now turn to pairs of points. Their contributions separate
naturally into pairs within a single cell and pairs belonging
to different cells. The former depend on the second-order
structure of the motifs, while the latter involve the joint
distribution of the marks. Recall that $n(m)$ is the expected number of points in a motif
of type $m$. Define its second moment by
\[
s(m)
:=
\int_{\mathcal N_f(Q)}
\nu(Q)^2\,P_m(\mathrm d\nu),
\qquad m\in\mathbb M.
\]
The map $m\mapsto s(m)$ is again measurable, since
$\nu\mapsto\nu(Q)^2$ is measurable and $m\mapsto P_m$
is a probability kernel.
Throughout the second-order analysis, we assume that
\begin{equation}\label{EQ_second_moment_integrability}
\overline s
:=
\int_{\mathbb M}s(m)\,\pi(\mathrm dm)
<
\infty.
\end{equation}
By Jensen's inequality, $n(m)^2\leq s(m)$, and hence
\[
\overline n
\leq
\left(
\int_{\mathbb M}n(m)^2\,\pi(\mathrm dm)
\right)^{1/2}
\leq
\overline s^{1/2}.
\]
In particular, $\widetilde\Phi$ has finite intensity.
The same assumption also guarantees the local
square-integrability needed to define its second-order measures.

\begin{lemm}[Local square-integrability]
Under \eqref{EQ_second_moment_integrability},
$\widetilde\Phi$ is locally square-integrable, that is,
\[
\E[\widetilde\Phi(B)^2]<\infty
\]
for every bounded Borel set $B\subseteq\mathbb R^d$.
\end{lemm}

\begin{proo}
For a bounded Borel set $B$, let
\[
K_B
:=
\left\{
k\in\mathbb Z^d:
(k+Q+Q)\cap B\neq\varnothing
\right\}.
\]
This set is finite and contains all indices of cells that
can contribute to $\theta_u\Phi(B)$, uniformly over $u\in Q$.
Consequently,
\[
\widetilde\Phi(B)
\leq
\sum_{k\in K_B}\Xi_k(Q)
\qquad\text{a.s}.
\]
The conditional motif law and stationarity of the mark field give
\[
\E[\Xi_k(Q)^2]
=
\E[s(M_k)]
=
\overline s.
\]
Thus, by the Cauchy--Schwarz inequality,
\[
\E[\widetilde\Phi(B)^2]
\leq
|K_B|\sum_{k\in K_B}\E[\Xi_k(Q)^2]
=
|K_B|^2\overline s
<
\infty.
\]
\qed
\end{proo}

\paragraph{Second-order measures of the motifs.}

For $m\in\mathbb M$, let $M_m^{(2)}$ and $M_{m,!}^{(2)}$
denote the second-order moment and factorial moment measures
of a motif with law $P_m$. They are determined by
\begin{align*}
M_m^{(2)}(C\times D)
&=
\int_{\mathcal N_f(Q)}
\nu(C)\nu(D)\,P_m(\mathrm d\nu),\\
M_{m,!}^{(2)}(C\times D)
&=
\int_{\mathcal N_f(Q)}
\bigl[\nu(C)\nu(D)-\nu(C\cap D)\bigr]
\,P_m(\mathrm d\nu),
\end{align*}
for $C,D\in\mathcal B(Q)$.
In particular, $M_m^{(2)}(Q\times Q)=s(m)$.
Under \eqref{EQ_second_moment_integrability}, these measures are finite
for $\pi$-almost every $m$, and their averages
\[
\overline M^{(2)}
:=
\int_{\mathbb M}M_m^{(2)}\,\pi(\mathrm dm),
\qquad
\overline M_!^{(2)}
:=
\int_{\mathbb M}M_{m,!}^{(2)}\,\pi(\mathrm dm)
\]
are finite measures satisfying
\[
\overline M^{(2)}(Q\times Q)=\overline s,
\qquad
\overline M_!^{(2)}(Q\times Q)=\overline s-\overline n.
\]
With $\Delta_Q(x):=(x,x)$, they are related by
\[
\overline M^{(2)}
=
\overline M_!^{(2)}+(\Delta_Q)_\#\overline\mu.
\]

\paragraph{Pairs within and between cells.}

The conditional sampling rule gives, for $C,D\in\mathcal B(Q)$,
\begin{equation}\label{EQ_conditional_second_moments_CIM}
\E\bigl[\Xi_k(C)\Xi_l(D)\mid\mathbf M\bigr]
=
\begin{cases}
M_{M_k}^{(2)}(C\times D), & k=l,\\[0.3em]
\mu_{M_k}(C)\mu_{M_l}(D), & k\neq l,
\end{cases}
\qquad\text{almost surely}.
\end{equation}
Thus, pairs within one cell involve the motif's second-order
measure, whereas conditional independence reduces the
contribution from distinct cells to a product of intensities.
For factorial moments within a cell,
\[
\E\bigl[
\Xi_k(C)\Xi_k(D)-\Xi_k(C\cap D)
\mid\mathbf M
\bigr]
=
M_{M_k,!}^{(2)}(C\times D)
\qquad\text{almost surely}.
\]

To express these contributions after translation, write
\[
C_{k,u}:=(C-k-u)\cap Q,
\qquad k\in\mathbb Z^d,\quad u\in Q,
\]
for a Borel set $C\subseteq\mathbb R^d$.
This is the portion of the reference cell mapped into $C$
by translation through $k+u$.

\begin{prop}[Second-order moment measures of the decorated lattice]
Assume \eqref{EQ_second_moment_integrability}.
For all bounded Borel sets $C,D\subseteq\mathbb R^d$,
\begin{align}
\alpha_{\widetilde\Phi}^{(2)}(C\times D)
={}&
\int_Q
\sum_{k\in\mathbb Z^d}
\overline M^{(2)}
\bigl(C_{k,u}\times D_{k,u}\bigr)
\,\mathrm du
\notag\\
&+
\int_Q
\sum_{\substack{k,l\in\mathbb Z^d\\k\neq l}}
\int_{\mathbb M^2}
\mu_m(C_{k,u})\mu_n(D_{l,u})
\,\pi_{l-k}(\mathrm dm,\mathrm dn)\,\mathrm du.
\label{EQ_second_moment_CIM_stationary}
\end{align}
Its second-order factorial moment measure is
\begin{align}
\alpha_{\widetilde\Phi,!}^{(2)}(C\times D)
={}&
\int_Q
\sum_{k\in\mathbb Z^d}
\overline M_!^{(2)}
\bigl(C_{k,u}\times D_{k,u}\bigr)
\,\mathrm du
\notag\\
&+
\int_Q
\sum_{\substack{k,l\in\mathbb Z^d\\k\neq l}}
\int_{\mathbb M^2}
\mu_m(C_{k,u})\mu_n(D_{l,u})
\,\pi_{l-k}(\mathrm dm,\mathrm dn)\,\mathrm du.
\label{EQ_second_factorial_moment_CIM_stationary}
\end{align}
\end{prop}

\begin{proo}
Independence of $U$ from the motif field and Tonelli's theorem
yield
\[
\alpha_{\widetilde\Phi}^{(2)}(C\times D)
=
\int_Q
\sum_{k,l\in\mathbb Z^d}
\E\bigl[\Xi_k(C_{k,u})\Xi_l(D_{l,u})\bigr]
\,\mathrm du.
\]
Since $C,D$ are bounded, only finitely many indices can
contribute, uniformly over $u\in Q$.

For $k=l$, taking expectations in
\eqref{EQ_conditional_second_moments_CIM} gives
\[
\E\bigl[\Xi_k(C_{k,u})\Xi_k(D_{k,u})\bigr]
=
\overline M^{(2)}
\bigl(C_{k,u}\times D_{k,u}\bigr).
\]
For $k\neq l$, conditional independence gives
\begin{align*}
\E\bigl[\Xi_k(C_{k,u})\Xi_l(D_{l,u})\bigr]
&=
\E\bigl[
\mu_{M_k}(C_{k,u})\mu_{M_l}(D_{l,u})
\bigr]=
\int_{\mathbb M^2}
\mu_m(C_{k,u})\mu_n(D_{l,u})
\,\pi_{l-k}(\mathrm dm,\mathrm dn),
\end{align*}
where we used stationarity of the mark field:
\[
\mathcal L(M_k,M_l)=\mathcal L(M_0,M_{l-k})=\pi_{l-k}.
\]
Splitting the sum into these two cases proves
\eqref{EQ_second_moment_CIM_stationary}.

For the factorial moment measure, subtract
$\E[\widetilde\Phi(C\cap D)]$.
Since
\[
(C\cap D)_{k,u}=C_{k,u}\cap D_{k,u},
\]
this replaces each same-cell contribution by
\begin{align*}
&\overline M^{(2)}(C_{k,u}\times D_{k,u})
-
\overline\mu(C_{k,u}\cap D_{k,u})=
\overline M_!^{(2)}(C_{k,u}\times D_{k,u}).
\end{align*}
The distinct-cell contributions are unchanged, since
self-pairs occur only within a single cell.
This proves
\eqref{EQ_second_factorial_moment_CIM_stationary}.
\qed
\end{proo}

\subsection{Bartlett's spectral measure}

We first compute the reduced second-order moment measures of
$\widetilde\Phi$. Passing to relative displacements removes
the common uniform shift and reduces the calculation to
pairs of motifs in cells indexed by $0$ and $r$.

\begin{prop}[Reduced second-order moment measures]
\label{P_reduced_second_moment}
Assume \eqref{EQ_second_moment_integrability}.
For every bounded Borel set $B\subseteq\mathbb R^d$,
\begin{align}
\alpha_{\widetilde\Phi,\mathrm{red}}^{(2)}(B)
={}&
\int_{Q\times Q}
\mathbf 1_B(x-y)\,
\overline M^{(2)}(\mathrm dx,\mathrm dy)
\notag\\
&+
\sum_{r\in\mathbb Z^d\setminus\{0\}}
\int_{\mathbb M^2}
\int_{Q\times Q}
\mathbf 1_B(x-y-r)\,
\mu_m(\mathrm dx)\mu_n(\mathrm dy)
\,\pi_r(\mathrm dm,\mathrm dn).
\label{EQ_reduced_second_moment_CIM}
\end{align}
Similarly,
\begin{align}
\alpha_{\widetilde\Phi,\mathrm{red},!}^{(2)}(B)
={}&
\int_{Q\times Q}
\mathbf 1_B(x-y)\,
\overline M_!^{(2)}(\mathrm dx,\mathrm dy)
\notag\\
&+
\sum_{r\in\mathbb Z^d\setminus\{0\}}
\int_{\mathbb M^2}
\int_{Q\times Q}
\mathbf 1_B(x-y-r)\,
\mu_m(\mathrm dx)\mu_n(\mathrm dy)
\,\pi_r(\mathrm dm,\mathrm dn).
\label{EQ_reduced_second_factorial_moment_CIM}
\end{align}
For bounded $B$, both sums have only finitely many
nonzero terms.
\end{prop}

\begin{proo}
By \eqref{EQ_reduced_second_moment}, independence of $U$
from the motif field, and Tonelli's theorem,
\begin{align*}
\alpha_{\widetilde\Phi,\mathrm{red}}^{(2)}(B)
&=
\int_Q
\sum_{k,l\in\mathbb Z^d}
\E\left[
\int_{Q\times Q}
\mathbf 1_Q(k+x+u)
\mathbf 1_B(k-l+x-y)
\,\Xi_k(\mathrm dx)\Xi_l(\mathrm dy)
\right]
\,\mathrm du.
\end{align*}
Set $r:=l-k$.
Stationarity of the motif field allows us to replace
$(\Xi_k,\Xi_{k+r})$ by $(\Xi_0,\Xi_r)$ inside the expectation.
For every fixed $x\in Q$,
\[
\sum_{k\in\mathbb Z^d}
\int_Q\mathbf 1_Q(k+x+u)\,\mathrm du
=
\sum_{k\in\mathbb Z^d}
\int_{k+x+Q}\mathbf 1_Q(z)\,\mathrm dz
=
\lambda_d(Q)
=
1.
\]
Consequently,
\begin{equation}\label{EQ_reduced_second_moment_general}
\alpha_{\widetilde\Phi,\mathrm{red}}^{(2)}(B)
=
\sum_{r\in\mathbb Z^d}
\E\left[
\int_{Q\times Q}
\mathbf 1_B(x-y-r)\,
\Xi_0(\mathrm dx)\Xi_r(\mathrm dy)
\right].
\end{equation}

For $r=0$, the expected product measure is
$\overline M^{(2)}$.
For $r\neq0$, conditional independence gives
\[
\E[\Xi_0\otimes\Xi_r]
=
\int_{\mathbb M^2}
\mu_m\otimes\mu_n\,\pi_r(\mathrm dm,\mathrm dn),
\]
where the equality is understood as an identity of measures.
Substitution into
\eqref{EQ_reduced_second_moment_general} proves
\eqref{EQ_reduced_second_moment_CIM}.

For the factorial measure, subtract
$\overline n\,\delta_0$.
Since
\[
\overline M^{(2)}-\overline M_!^{(2)}
=
(\Delta_Q)_\#\overline\mu,
\]
the same-cell term changes by
\[
\int_{Q\times Q}
\mathbf 1_B(x-y)\,
(\Delta_Q)_\#\overline\mu(\mathrm dx,\mathrm dy)
=
\overline n\,\mathbf 1_B(0).
\]
The distinct-cell terms are unchanged, proving
\eqref{EQ_reduced_second_factorial_moment_CIM}.

Finally, $x-y\in Q-Q=(-1,1)^d$, so
$x-y-r\in B$ requires
\[
r\in(Q-Q-B)\cap\mathbb Z^d.
\]
This set is finite whenever $B$ is bounded.
\qed
\end{proo}

The following theorem gives an explicit spectral description of
the decorated lattice. It separates the fluctuations within
individual motifs, the dependence transmitted by the mark field,
and the periodic contribution of the averaged motif.
This decomposition will allow us to study large-scale
fluctuations through the behaviour of the spectrum near the origin.

For $m\in\mathbb M$, define
\[
\widehat\mu_m(k)
:=
\int_Q e^{-ik\cdot x}\,\mu_m(\mathrm dx),
\qquad k\in\mathbb R^d,
\]
and
\[
\widehat M_m^{(2)}(k)
:=
\int_{Q\times Q}
e^{-ik\cdot(x-y)}\,M_m^{(2)}(\mathrm dx,\mathrm dy).
\]
These quantities are finite for almost all $m \in \mathbb{M}$ by
\eqref{EQ_second_moment_integrability}.
The latter is the Fourier transform in the displacement
variable $x-y$.
Further, let \[
\widehat{\overline\mu}(k)
=
\int_{\mathbb M}\widehat\mu_m(k)\,\pi(\mathrm dm).
\]

For $r\in\mathbb Z^d$, let
\[
\eta_r:=\pi_r-\pi\otimes\pi
\]
measure the deviation of the two-mark distribution from
independence, and define
\begin{equation}\label{EQ_Gamma_CIM}
\Gamma_r(k)
:=
\C\bigl(\widehat\mu_{M_0}(k),\widehat\mu_{M_r}(k)\bigr)
=
\int_{\mathbb M^2}
\widehat\mu_a(k)\overline{\widehat\mu_b(k)}
\,\eta_r(\mathrm da,\mathrm db).
\end{equation}
Finally, write
\[
Z_j(k)
:=
\int_Q e^{-ik\cdot x}\,\Xi_j(\mathrm dx),
\qquad j\in\mathbb Z^d,
\]
for the Fourier amplitude of the motif in cell $j$.

\begin{theo}[Bartlett spectral measure of the decorated lattice]
\label{L_Bartlett_CIM}
Assume \eqref{EQ_second_moment_integrability} and
\[
\sum_{r\in\mathbb Z^d}
\int_{\mathbb M^2}
n(a)n(b)\,|\eta_r|(\mathrm da,\mathrm db)
<
\infty,
\]
where $|\eta_r|$ denotes the total variation measure of $\eta_r$.
Then
\begin{equation}\label{EQ_Bartlett_CIM}
\widehat\beta_{\widetilde\Phi}(\mathrm dk)
=
S(k)\,\mathrm dk
+
(2\pi)^d
\sum_{z\in\mathbb Z^d\setminus\{0\}}
\left|\widehat{\overline\mu}(2\pi z)\right|^2
\,\delta_{2\pi z}(\mathrm dk),
\end{equation}
where
\begin{align}
S(k)
={}&
\int_{\mathbb M}
\left(
\widehat M_m^{(2)}(k)
-
|\widehat\mu_m(k)|^2
\right)\pi(\mathrm dm)
+
\sum_{r\in\mathbb Z^d}
e^{ik\cdot r}\Gamma_r(k).
\label{EQ_S_CIM}
\end{align}
Equivalently,
\begin{equation}\label{EQ_S_CIM_covariance_representation}
S(k)
=
\sum_{r\in\mathbb Z^d}
e^{ik\cdot r}\,
\C\bigl(Z_0(k),Z_r(k)\bigr).
\end{equation}
Both series converge absolutely and uniformly in
$k\in\mathbb R^d$.
The function $S$ is bounded, continuous, even, and nonnegative.
\end{theo}

The first term in \eqref{EQ_S_CIM} is the expected conditional
variance of a cell's Fourier amplitude. The second describes
the correlations of the mark-selected mean motifs, including
their variation at displacement zero.
The Bragg peaks depend only on the averaged motif intensity
$\overline\mu$.
Thus, under the summability assumption, dependence contributes
to the continuous part of the spectrum, while the remaining
periodic component is determined by the average motif.

\begin{proo}
For a finite signed measure $K$ on $Q\times Q$ and
$r\in\mathbb Z^d$, let $\mathcal T_rK$ be its pushforward
under $(x,y)\mapsto x-y-r$, that is,
\[
\mathcal T_rK(B)
=
\int_{Q\times Q}\mathbf 1_B(x-y-r)\,
K(\mathrm dx,\mathrm dy).
\]
All integrals over marks below are understood up to
$\pi$-null sets. Such exceptional sets are also null for
each marginal of $\pi_r$ and $|\eta_r|$.

By Proposition~\ref{P_reduced_second_moment},
\begin{align*}
\alpha_{\widetilde\Phi,\mathrm{red}}^{(2)}
={}&
\mathcal T_0\overline M^{(2)}+
\sum_{r\in\mathbb Z^d\setminus\{0\}}
\int_{\mathbb M^2}
\mathcal T_r(\mu_a\otimes\mu_b)
\,\pi_r(\mathrm da,\mathrm db).
\end{align*}
We decompose this measure as
\[
\alpha_{\widetilde\Phi,\mathrm{red}}^{(2)}
=
A+B+C,
\]
where
\begin{align*}
A
&:=
\sum_{r\in\mathbb Z^d}
\mathcal T_r(\overline\mu\otimes\overline\mu),\\
B
&:=
\int_{\mathbb M}
\mathcal T_0\bigl(M_m^{(2)}-\mu_m\otimes\mu_m\bigr)
\,\pi(\mathrm dm),\\
C
&:=
\sum_{r\in\mathbb Z^d}
\int_{\mathbb M^2}
\mathcal T_r(\mu_a\otimes\mu_b)
\,\eta_r(\mathrm da,\mathrm db).
\end{align*}
Indeed, $\pi_r=\pi\otimes\pi+\eta_r$, while
$\pi_0=\mathcal L(M_0,M_0)$ is concentrated on the diagonal.
Thus, the $r=0$ contribution in $A+C$ is
\[
\int_{\mathbb M}\mathcal T_0(\mu_m\otimes\mu_m)
\,\pi(\mathrm dm),
\]
and adding $B$ replaces it by
$\mathcal T_0\overline M^{(2)}$.

The measure $A$ represents the periodic background, $B$ the
conditional fluctuations within motifs, and $C$ the
correlations of their conditional means.

\medskip
\noindent\emph{The periodic contribution.}
The measure $A$ is locally finite and $\mathbb Z^d$-periodic.
Each summand is concentrated on $Q-Q-r$ and has mass
$\overline n^2$. A translate of any fixed bounded set
intersects a uniformly bounded number of these sets.
Hence $A$ is translation bounded and defines a tempered
distribution.

Since $\overline\mu$ is finite and concentrated on the bounded
set $Q$, all derivatives of $\widehat{\overline\mu}$ are bounded.
Therefore,
$h|\widehat{\overline\mu}|^2\in\mathcal S(\mathbb R^d)$
whenever $h\in\mathcal S(\mathbb R^d)$.
Poisson summation, with our Fourier convention, gives
\[
\sum_{r\in\mathbb Z^d}e^{-ik\cdot r}
=
(2\pi)^d\sum_{z\in\mathbb Z^d}\delta_{2\pi z}
\]
in $\mathcal S'(\mathbb R^d)$.
Consequently,
\begin{align*}
\widehat A(h)
&=
(2\pi)^d A(\check h)=
\left\langle
\sum_{r\in\mathbb Z^d}e^{-ik\cdot r},
h(k)|\widehat{\overline\mu}(k)|^2
\right\rangle=
(2\pi)^d
\sum_{z\in\mathbb Z^d}
h(2\pi z)|\widehat{\overline\mu}(2\pi z)|^2.
\end{align*}
Thus,
\[
\widehat A
=
(2\pi)^d
\sum_{z\in\mathbb Z^d}
|\widehat{\overline\mu}(2\pi z)|^2\,\delta_{2\pi z}.
\]

\medskip
\noindent\emph{The fluctuation contributions.}
The measures $B$ and $C$ have finite total variation.
Indeed,
\[
\|B\|_{\mathrm{TV}}
\leq
\int_{\mathbb M}\bigl(s(m)+n(m)^2\bigr)\,\pi(\mathrm dm)
\leq
2\overline s,
\]
and
\[
\|C\|_{\mathrm{TV}}
\leq
\sum_{r\in\mathbb Z^d}
\int_{\mathbb M^2}
n(a)n(b)\,|\eta_r|(\mathrm da,\mathrm db)
<
\infty.
\]
In particular, the defining series for $C$ converges in
total variation.

Both measures are symmetric under reflection.
For $B$, this follows from invariance of $M_m^{(2)}$ and
$\mu_m\otimes\mu_m$ under interchange of the coordinates.
For $C$, stationarity gives
\[
\eta_{-r}=\tau_\#\eta_r,
\qquad \tau(a,b):=(b,a),
\]
and reflection sends
$\mathcal T_r(\mu_a\otimes\mu_b)$ to
$\mathcal T_{-r}(\mu_b\otimes\mu_a)$.
Their Fourier transforms therefore have densities
$\mathcal F(B)$ and $\mathcal F(C)$, respectively.

Direct calculation yields
\[
\widehat B(\mathrm dk)
=
\left[
\int_{\mathbb M}
\left(
\widehat M_m^{(2)}(k)-|\widehat\mu_m(k)|^2
\right)\pi(\mathrm dm)
\right]\mathrm dk.
\]
Similarly,
\[
\mathcal F\bigl(\mathcal T_r(\mu_a\otimes\mu_b)\bigr)(k)
=
e^{ik\cdot r}
\widehat\mu_a(k)\overline{\widehat\mu_b(k)},
\]
so
\[
\widehat C(\mathrm dk)
=
\left[
\sum_{r\in\mathbb Z^d}e^{ik\cdot r}\Gamma_r(k)
\right]\mathrm dk.
\]
The series converges absolutely and uniformly, because
\[
\sup_{k\in\mathbb R^d}|\Gamma_r(k)|
\leq
\int_{\mathbb M^2}
n(a)n(b)\,|\eta_r|(\mathrm da,\mathrm db).
\]

Combining the three transforms gives
\[
\widehat\alpha_{\widetilde\Phi}(\mathrm dk)
=
S(k)\,\mathrm dk
+
(2\pi)^d
\sum_{z\in\mathbb Z^d}
|\widehat{\overline\mu}(2\pi z)|^2
\,\delta_{2\pi z}(\mathrm dk).
\]
Since $\widetilde\Phi$ has intensity $\overline n$,
\[
\widehat\beta_{\widetilde\Phi}
=
\widehat\alpha_{\widetilde\Phi}
-
(2\pi)^d\overline n^2\delta_0.
\]
Moreover, $\widehat{\overline\mu}(0)=\overline n$, so this
subtraction cancels precisely the $z=0$ term.
This proves \eqref{EQ_Bartlett_CIM}.

The density $S=\mathcal F(B)+\mathcal F(C)$ is bounded and
continuous, since $B,C$ are finite signed measures.
Their symmetry implies that $S$ is real-valued and even.
Positivity of the Bartlett spectral measure implies
$S\geq0$ almost everywhere, and continuity then gives
$S(k)\geq0$ for every $k$.

\medskip
\noindent\emph{The covariance representation.}
The conditional motif law gives
\[
\E[Z_j(k)\mid\mathbf M]=\widehat\mu_{M_j}(k),
\qquad
\E[|Z_j(k)|^2\mid\mathbf M]
=
\widehat M_{M_j}^{(2)}(k).
\]
For $r\neq0$, conditional independence and the law of total
covariance yield
\[
\C\bigl(Z_0(k),Z_r(k)\bigr)=\Gamma_r(k).
\]
For $r=0$, the law of total variance gives
\[
\V(Z_0(k))
=
\int_{\mathbb M}
\left(
\widehat M_m^{(2)}(k)-|\widehat\mu_m(k)|^2
\right)\pi(\mathrm dm)
+
\Gamma_0(k).
\]
These identities prove
\eqref{EQ_S_CIM_covariance_representation}.
Its absolute and uniform convergence follows from the bound
on $\Gamma_r$ above and
\[
\V(Z_0(k))
\leq
\E[|Z_0(k)|^2]
\leq
\E[\Xi_0(Q)^2]
=
\overline s.
\]
\qed
\end{proo}

\subsection{Hyperuniformity and \textit{p}-uniformity of decorated lattices}
Let
\[
N_j:=\Xi_j(Q),
\qquad j\in\mathbb Z^d,
\]
be the number of points placed in the cell $j+Q$. Since
\[
Z_j(0)=N_j,
\]
the covariance representation
\eqref{EQ_S_CIM_covariance_representation} gives
\begin{equation}
    \label{e:S-0}
S(0)
=
\sum_{r\in\mathbb Z^d}
\C(N_0,N_r).
\end{equation}
The series is absolutely convergent under the assumptions of
Theorem~\ref{L_Bartlett_CIM}. This gives us direct conditions for hyperuniformity of decorated lattices. 
Further, it turns out that this formula relates the hyperuniformity of the decorated lattice to the hyperuniformity of a random field on the lattice. Hence, we extend the concept of hyperuniformity and $p$-uniformity to random fields at this point.

\paragraph{Hyperuniformity of stationary random fields.}

Suppose that $Y=(Y_j)_{j\in\mathbb Z^d}$ is a  real-valued stationary
random field with $\E[Y_0^2]<\infty$.
By the Herglotz theorem, its covariance function admits a
unique finite positive spectral measure $\widehat\beta_Y$
on the torus
\(
\mathbb T^d:=\mathbb R^d/(2\pi\mathbb Z^d).
\)
With our Fourier convention, this measure satisfies
\[
\C\left(
\sum_{j\in\mathbb Z^d}f(j)Y_j,
\sum_{j\in\mathbb Z^d}g(j)Y_j
\right)
=
\frac{1}{(2\pi)^d}
\int_{\mathbb T^d}
\widehat f(k)\overline{\widehat g(k)}
\,\widehat\beta_Y(\mathrm dk)
\]
for finitely supported functions
$f,g:\mathbb Z^d\to\mathbb C$, where
\[
\widehat f(k)
:=
\sum_{j\in\mathbb Z^d}e^{-ik\cdot j}f(j).
\]
We identify $\mathbb T^d$ with $[-\pi,\pi)^d$, so that
$B_\varepsilon$, for $0<\varepsilon<\pi$, is an ordinary
Euclidean ball around the origin.
\begin{defn}[$p$-uniform random fields]
Let $Y=(Y_j)_{j \in \Z^d}$ be a real-valued stationary random field with $\E[Y^2_0]<\infty$. In analogy with Definition~\ref{d:p-uniformity}, we call
$Y$ \emph{$p$-uniform}, for $p\in[-d,\infty)$, if
\[
\frac{\widehat\beta_Y(B_\varepsilon)}
{\varepsilon^d}
=
O(\varepsilon^p)
\qquad\text{as }\varepsilon\downarrow0,
\]
and \emph{beyond $p$-uniform} if the right-hand side is
$o(\varepsilon^p)$. We call $Y$ \emph{hyperuniform} if
it is beyond $0$-uniform.    
\end{defn}
 If the covariance function is absolutely summable, then
$\widehat\beta_Y$ has the continuous, nonnegative density
\[
S_Y(k)
=
\sum_{r\in\mathbb Z^d}
e^{ik\cdot r}\C(Y_0,Y_r).
\]
Consequently,
\[
Y\text{ is hyperuniform}
\quad\Longleftrightarrow\quad
S_Y(0)
=
\sum_{r\in\mathbb Z^d}\C(Y_0,Y_r)
=
0.
\]

Under the assumptions of Theorem~\ref{L_Bartlett_CIM},
the count field $N=(N_j)_{j\in\mathbb Z^d}$ has absolutely
summable covariances, and \eqref{e:S-0} gives
\(
S(0)=S_N(0).
\)
Moreover, the Bartlett spectral measure of $\widetilde\Phi$
has continuous density $S$ in a neighborhood of the origin.
Consequently,
\[
\widetilde\Phi\text{ is hyperuniform}
\quad\Longleftrightarrow\quad
N\text{ is hyperuniform}.
\]
The following proposition expresses this criterion in terms
of the motif laws and the mark field.

\begin{prop}[Hyperuniformity of decorated lattices]
\label{P_hyperuniformity_decorated}
Under the assumptions of Theorem~\ref{L_Bartlett_CIM}, hyperuniformity of $\widetilde\Phi$ is equivalent to $S(0)=0$ with 
\begin{equation}\label{EQ_S_CIM_zero_decomposition}
S(0)
=
\E[\V(N_0\mid M_0)]
+
\sum_{r\in\mathbb Z^d}
\C\left(
\mathbb E[N_0\mid M_0],
\mathbb E[N_r\mid M_r]
\right).
\end{equation}

Consequently, $\widetilde\Phi$ is hyperuniform if and only if
\(
\V(N_0\mid M_0)=0 \text{ almost surely}
\)
and the field $(\E[N_j\mid M_j])_{j\in\mathbb Z^d}$ is hyperuniform. In particular,

\begin{enumerate}
\item
If $\E[\V(N_0\mid M_0)]>0$, then $\widetilde\Phi$
is not hyperuniform.

\item If the number of points is determined by the mark, that is,
\[
N_j=c(M_j),
\]
almost surely for a measurable
$c:\mathbb M\to\mathbb N_0$
then
\[
S(0)
=
\sum_{r\in\mathbb Z^d}
\C\bigl(c(M_0),c(M_r)\bigr).
\]
Hence the decorated lattice is hyperuniform if and only if the
mark-induced count field $(c(M_j))_{j\in\mathbb Z^d}$ is hyperuniform

\item Specifically, if every motif has the same deterministic number of points,
\[
N_j=c
\qquad\text{a.s.},
\]
then
\(
S(0)=0,
\)
and the decorated lattice is hyperuniform.

\item If the mark field is i.i.d., then
\[
S(0)
=
\V(N_0).
\]
Thus, in the i.i.d.\ case, hyperuniformity holds if and only if the
cell count $N_0$ is deterministic.
\end{enumerate}
\end{prop}

\begin{proo}
By Theorem~\ref{L_Bartlett_CIM}, $S$ is continuous and
\[
\widehat\beta_{\widetilde\Phi}(B_\varepsilon)
=
\int_{B_\varepsilon}S(k)\,\mathrm dk,
\qquad 0<\varepsilon<2\pi.
\]
Consequently,
\[
\lim_{\varepsilon\to 0}
\frac{\widehat\beta_{\widetilde\Phi}(B_\varepsilon)}
     {\varepsilon^d}
=
\lambda_d(B_1)S(0),
\]
so $\widetilde\Phi$ is hyperuniform if and only if $S(0)=0$. 

For $r\neq0$, conditional independence of the motifs and
the law of total covariance therefore imply
\[
\C(N_0,N_r)
=
\C\left(
\E[N_0\mid M_0],\E[N_r\mid M_r]
\right).
\]
For $r=0$, the law of total variance gives
\[
\V(N_0)
=
\E[\V(N_0\mid M_0)]
+
\V\bigl(\E[N_0\mid M_0]\bigr).
\]
Substituting these identities into \eqref{e:S-0}
proves \eqref{EQ_S_CIM_zero_decomposition}.

The field $(\E[N_j\mid M_j])_{j\in\mathbb Z^d}$ is
stationary and square-integrable. Its covariances are
absolutely summable: at nonzero displacements they equal
the count covariances, and
\[
\V\bigl(\E[N_0\mid M_0]\bigr)\leq\V(N_0)<\infty.
\]
It therefore has a continuous, nonnegative spectral
density whose value at zero is
\[
\sum_{r\in\mathbb Z^d}
\C\left(
\E[N_0\mid M_0],\E[N_r\mid M_r]
\right).
\]
This sum is nonnegative and vanishes exactly when the
conditional-mean field is hyperuniform.

Both terms in \eqref{EQ_S_CIM_zero_decomposition} are
thus nonnegative. The first vanishes exactly when
$\V(N_0\mid M_0)=0$ almost surely, and the second
vanishes exactly when
$(\E[N_j\mid M_j])_{j\in\mathbb Z^d}$ is hyperuniform.
This proves the stated characterization and the first
three consequences.

Finally, if the mark field is i.i.d., its conditional-mean
field is also i.i.d., so all terms with $r\neq0$ vanish.
The law of total variance then yields
\[
S(0)
=
\E[\V(N_0\mid M_0)]
+
\V\bigl(\E[N_0\mid M_0]\bigr)
=
\V(N_0).
\]
Hence hyperuniformity holds exactly when $N_0$ is
deterministic.
\qed
\end{proo}

The decomposition shows that residual randomness in the
number of points within a cell cannot be cancelled by
negative correlations between marks.
Hyperuniformity therefore requires the cell count to be
determined by its mark, together with a balance of the
remaining fluctuations across cells.
Under the assumptions above, the positions of points within
the motifs do not affect whether the process is hyperuniform.
They can, however, affect its finer order of $p$-uniformity.

\begin{prop}[$p$-uniformity of decorated lattices]
Assume the hypotheses of Theorem~\ref{L_Bartlett_CIM}.
Then, for
$p\in[-d,\infty)$, $\widetilde\Phi$ is $p$-uniform
if 
\[
S(k)
=
O(| k|^p)
\qquad\text{as }|k|\to 0,
\]
It is beyond $p$-uniform if the right-hand side is
$o(|k|^p)$.
For every $p\in[-d,0)$, the process is automatically
beyond $p$-uniform.
\end{prop}

\begin{proo}
For $0<\varepsilon<2\pi$, identity
\eqref{EQ_Bartlett_CIM} gives
\[
\frac{\widehat\beta(B_\varepsilon)}
{\varepsilon^d}
=
\frac{1}{\varepsilon^d}
\int_{B_\varepsilon}
S(k)\,\mathrm dk.
\]
First let $p\geq0$. If $S(k)=O(\lVert k\rVert^p)$, then
\[
\int_{B_\varepsilon}S(k)\,\mathrm dk
=
O(\varepsilon^{d+p}),
\]
and the implication follows. For $p\in[-d,0)$, boundedness of $S$ gives
\[
\frac{\widehat\beta_{\widetilde\Phi}(B_\varepsilon)}
{\varepsilon^d}
=
O(1)
=
o(\varepsilon^p).
\]
Thus, the process is automatically beyond $p$-uniform
throughout this range.

\qed
\end{proo}

\subsection{Extension to general lattices}
The restriction to the cubic lattice is inessential.
Let $A$ be an invertible $d\times d$ matrix and put
$v:=|\det A|$. The image process
$\widetilde\Phi_A:=A_\#\widetilde\Phi$ is a stationary
decoration of a lattice $A \Z^d$ with cells $Ak+AQ$, $k\in\mathbb Z^d$,
and intensity $\overline n/v$. Its Bartlett spectral
measure satisfies
\[
\widehat\beta_{\widetilde\Phi_A}
=
v^{-2}(A^{-\mathsf T})_\#
\widehat\beta_{\widetilde\Phi}.
\]
In particular, its continuous spectral density is
$S_A(k)=v^{-1}S(A^{\mathsf T}k)$, and its Bragg peaks
lie on $2\pi A^{-\mathsf T}\mathbb Z^d\setminus\{0\}$.
Since invertible linear maps send balls to comparable
ellipsoids, $p$-uniformity and beyond $p$-uniformity,
including hyperuniformity, are preserved.

\section{Central limit theorems for linear statistics}\label{sec:5}

Given that the generating random field is sufficiently mixing, the decorated lattices follow a central limit theorem. Here, we work with the following notion of an $\alpha$-mixing coefficient. Let $(M_k)_{k\in \Z^d}$ be the mark field, then 
\begin{equation}
    \alpha_{n,l}(r) := \sup\{\alpha(\sigma((M_k)_{k\in K_1}), \sigma((M_k)_{k\in K_2})): |K_1|\leq n, |K_2|\leq l, d(K_1, K_2) \geq r\},
\end{equation}
for $n,l\in\N\cup\{\infty\}, r>0$, where 
\begin{align*}
    \alpha(\mathcal A, \mathcal B) &= \sup\{|\mathbb P(A\cap B)-\mathbb P(A)\mathbb P(B)|: A\in\mathcal A, B\in\mathcal B\},\\
    d(K_1, K_2) &= \inf \{\|k_1 - k_2\|:k_1\in K_1, k_2\in K_2\},
\end{align*}
for sub-$\sigma$-algebras $\mathcal A, \mathcal B$ and $K_1,K_2\subseteq \Z^d$. The particular mixing condition we assume is that there are $\varepsilon>0, c>0$ such that
\begin{equation}\label{e:mixing-condition}
    \max(\alpha_{2,2}(r), \alpha_{1,\infty}(r))\leq c r^{-(d+\varepsilon)},\quad r>0.
\end{equation}
Now, let $\tilde\Phi = \sum_{k\in\Z^d} \theta_{k+U}\Xi_k, U\sim \mathcal U([0,1)^d),$ be a stationarized decorated lattice with mark field $(M_k)_{k\in \Z^d}$. Since $\Xi_0$ and $\Xi_k$ are independent for $k\in\Z^d\setminus\{0\}$ given their respective marks $M_0$ and $M_k$, the alpha mixing coefficient of the random field $(\Xi_k)_{k\in\Z^d}$ is bounded by that of the mark field $(M_k)_{k\in\Z^d}$. Let $(N_k)_{k\in\Z^d}$ with $N_k:=\Xi_k([0,1)^d), k\in\Z^d,$ be the random field that counts the number of points in each cell. We assume that 
\begin{equation}\label{e:moment-assumption}
    \E[N_0^m] < \infty
\end{equation}
for some $m>2+2\frac{d}{\varepsilon}$ for the $\varepsilon$ from \eqref{e:mixing-condition}. Then a standard inequality between the covariance and the $\alpha$-mixing coefficient yields 
\begin{equation}\label{e:covariance-integrability}
    \sum_{k\in\Z^d} \big|\C\big(\Xi_0(f), \Xi_k(f)\big)\big| \leq 8 \|f\|_\infty^2 \E[N_0^{m}]^{\frac{2}{m}}\sum_{k\in\Z^d} \alpha_{1,1}(\|k\|)^{1-\frac{2}{m}} < \infty
\end{equation}
for any bounded function $f:[0,1)^d\to\R$; see~\cite{Do94}. Hence, the random field $(N_k)_{k\in\Z^d}$ is a locally square-integrable random field with absolutely summable covariance function. In the non-hyperuniform regime, the asymptotic behavior of $\tilde\Phi$ is dominated by $(N_k)_{k\in\Z^d}$ and we obtain the following central limit theorem.
\begin{proposition}\label{p:classical-clt}
    Suppose that $\tilde\Phi$ is a stationarized decorated lattice with a mark field that satisfies the $\alpha$-mixing condition \eqref{e:mixing-condition} for some $\varepsilon>0$. Let $(N_k)_{k\in\Z^d}$ be the random field that counts the number of points per cell, and assume that \eqref{e:moment-assumption} holds for some $m>2+2\frac{d}{\varepsilon}$. Then $\tilde\Phi$ is hyperuniform if and only if $(N_k)_{k\in\Z^d}$ is hyperuniform. Further, given $f\in L^1\cap L^2$ and constants $c,\vartheta>0$ such that
    \begin{equation}
        |\hat f(k)| \leq c (1 + \|k\|)^{-\frac{d+\vartheta}{2}}, \quad k\in\R^d,
    \end{equation}
    we obtain
    \begin{equation}
        \frac{1}{r^{\frac{d}{2}}} \bigg(\int f(x/r)\, \tilde\Phi(\mathrm dx)-r^d\lambda_d(f)\E[N_0]\bigg) \xrightarrow{\mathcal D} \mathcal N(0, \sigma^2),
    \end{equation}
    where
    \begin{equation}
        \sigma^2 = \|f\|_2^2\sum_{k\in\Z^d} \C(N_0, N_k).
    \end{equation}
    If $f\neq0$, $\sigma^2>0$ is equivalent to $(N_k)_{k\in\Z^d}$ not being hyperuniform.
\end{proposition}
The test function $f$ can be chosen to be any indicator function over a Fourier smooth set, e.g., $f:=\I_{B_1}$; see \cite{BH24} for further examples. The proof is a simpler version of that of the following Proposition \ref{p:clt}. We should also note that Proposition \ref{p:classical-clt} establishes, that the characterization of hyperuniformity of $\tilde\Phi$ in terms of hyperuniformity of $(N_k)_{k\in\Z^d}$ from Proposition \ref{P_hyperuniformity_decorated} also holds under these alternate mixing conditions. In the case that $\tilde\Phi$ is hyperuniform, the picture is more complicated. If we consider the canonical transport from $\Z^d+U$ to $\tilde\Phi$ and calculate the Taylor expansion that was developed in \cite{LK26}, we see that $(N_k)_{k\in\Z^d}$ is the zeroth term of that expansion. Given the mixing assumption \eqref{e:mixing-condition}, the remaining terms of the expansion can start to dominate the asymptotic behavior of $\tilde\Phi$ when $(N_k)_{k\in\Z^d}$ is beyond $2$-uniform. We assume this in the following theorem.  First, we apply the expansion and Slutsky's lemma to reduce the problem to deriving a central limit theorem for the first term of the expansion. Then, we use a version of Bolthausens central limit theorem that allows for a weighted sum of random variables; see~\cite{JP09}. For the specific case that $\tilde\Phi$ is an independently perturbed lattice, a similar technique was applied in~\cite{Ma25}.

It may be surprising that the central limit theorem holds even for $d\in\{1,2\}$, where the variances of the ``partial sums" converge to $0$ or are uniformly bounded.

\begin{proposition}\label{p:clt}
    Suppose that $\tilde\Phi$ is a stationarized decorated lattice with a mark field that satisfies the $\alpha$-mixing condition \eqref{e:mixing-condition} for some $\varepsilon>0$. Let $(N_k)_{k\in\Z^d}$ be the random field that counts the number of points per cell and assume that it satisfies \eqref{e:moment-assumption} for some $m>2+2\frac{d}{\varepsilon}$ and that it is beyond $2$-uniform. Then $\widetilde{\Phi}$ is $2$-uniform. Further, given $f\in L^1$ and constants $c,\vartheta>0$ such that
    \begin{align}\label{e:real-space-decay}
        |f(x)| &\leq c (1 + \|x\|)^{-(d+1+\vartheta)}, \quad x\in\R^d,\\\label{e:fourier-space-decay}
        |\hat f(k)| &\leq c (1 + \|k\|)^{-\frac{d+2+\vartheta}{2}}, \quad k\in\R^d,
    \end{align}
    we obtain
    \begin{equation}\label{e:clt}
        \frac{1}{r^{\frac{d}{2}-1}} \bigg(\int f(x/r)\, \tilde\Phi(\mathrm dx)-r^d\lambda_d(f)\E[N_0]\bigg) \xrightarrow{\mathcal D} \mathcal N(0, \sigma^2),
    \end{equation}
    where
    \begin{equation}\label{e:asymptotic-variance}
        \sigma^2 = \sum_{n,l=1}^d\bigg(\langle \partial_n f , \partial_l f\rangle \sum_{k\in\Z^d} \C\bigg(\int x_n\, \Xi_0(\mathrm dx), \int x_l\, \Xi_k(\mathrm dx)\bigg)\bigg).
    \end{equation}
    If $f\neq0$, then $\sigma^2>0$ is equivalent to
    \begin{equation*}
        \sum_{k\in\Z^d} \C\bigg(\int x_n\, \Xi_0(\mathrm dx), \int x_n\, \Xi_k(\mathrm dx)\bigg) >0\text{ for at least one } n\in\{1,...,d\}.
    \end{equation*}
\end{proposition}

\begin{remark}
    The conditions on $f$ do not allow for it to be a non-trivial indicator function. An admissible example is the Gauss kernel $f(x):=(2\pi)^{-\frac{d}{2}}e^{-\frac{x^T x}{2}}, x\in\R^d$. In this case, the formula \eqref{e:asymptotic-variance} simplifies to
\begin{equation}
    \sigma^2 = \frac{1}{2^{d+1}\pi^\frac{d}{2}}\sum_{n=1}^d\sum_{k\in\Z^d} \C\bigg(\int x_n\, \Xi_0(\mathrm dx), \int x_n\, \Xi_k(\mathrm dx)\bigg).
\end{equation}
\end{remark}

\begin{proo}[of Proposition \ref{p:clt}]
    The proof is divided into three steps. First, we derive the variance asymptotics for a more restricted class of test functions. Then we prove the central limit theorem in this setting. In the third part, we generalize the class of test functions.
    
    \smallskip\noindent
    \underline{Variance asymptotics:}
    \par\nopagebreak[4]
    \smallskip\noindent
    Let us consider the random transport $K$ that moves mass from each point in the stationarized lattice $\Gamma:=\sum_{k\in\Z^d}\delta_{k+U}$ to the mass of $\tilde\Phi$ within the corresponding cell, whereby
    \begin{align*}
        K(k+U, \cdot) &= \tilde\Phi|_{k+U+[0,1)^d},\\
        K\Gamma &= \tilde\Phi.
    \end{align*}
    See~\cite[Subsection 2.2]{LK26} for a proper introduction of transports. Note that the mass of $\tilde\Phi$ in each cell need not be $1$, so this transport is allowed to create and destroy mass. However, by the moment assumption \eqref{e:moment-assumption}, the mass created in the typical cell has a finite second moment, and further, the distance over which it moves is uniformly bounded. Let $f\in\mathcal C_c^\infty$. We can apply~\cite[Theorem 4.14]{LK26} to obtain an expansion for the asymptotic behavior of $\int f(x/r)\, \tilde\Phi(\mathrm dx)$ as $r\to\infty$. Since the mixing assumption \eqref{e:mixing-condition} implies \eqref{e:covariance-integrability}, we obtain that the random field $(\int x_1^{n_1}\cdots x_d^{n_d}\, \Xi_k(\mathrm dx))_{k\in\Z^d}$ has an absolutely summable covariance function for any $n_1,...,n_d\in\N_0$. Therefore, all but the zeroth and first term of the expansion drop out since all components of the tensor-valued random measure
    \begin{equation*}
        \Psi_q := \sum_{k\in\Z^d} \bigg(\int x^{\otimes q}\, \Xi_k(\mathrm dx) \bigg)\delta_{k+U}
    \end{equation*}
    are $0$-uniform for any $q\in\N_0$, where $x^{\otimes q}$ is the $q$-times tensor product of $x$ with itself; see~\cite[Proposition 3.12]{LK26}. We obtain
    \begin{align*}
        \int f(x/r)\, \tilde\Phi(\mathrm dx) &= \Psi_0(f(\cdot/r)) + \frac{1}{r} \Psi_1(f'(\cdot/r)) + o(r^{\frac{d}{2}-1})\\
        &=\sum_{k\in\Z^d} f((k+U)/r) N_k + \frac{1}{r}\sum_{n=1}^d\sum_{k\in\Z^d} \partial_n f((k+U)/r) \int x_n\, \Xi_k(\mathrm dx) + o(r^{\frac{d}{2}-1})
    \end{align*}
    in $L^2$ as $r\to\infty$. Further, since $(N_k)_{k\in\Z^d}$ is assumed to be beyond $2$-uniform, the variance of the zeroth term is $o(r^{d-2})$, and 
    \begin{equation}\label{e:leading-order-decomposition}
        \int f(x/r)\, \tilde\Phi(\mathrm dx) = r^d\lambda_d(f)\E[N_0]+\frac{1}{r}\sum_{n=1}^d\sum_{k\in\Z^d} \partial_n f((k+U)/r) \int x_n\, \Xi_k(\mathrm dx) + o(r^{\frac{d}{2}-1})
    \end{equation}
    in $L^2$ as $r\to\infty$. Then, by the dominated convergence theorem,
    \begin{align}
        &\lim_{r\to\infty}\V\bigg(\frac{1}{r^{\frac{d}{2}-1}}\int f(x/r)\, \tilde\Phi(\mathrm dx)\bigg) \nonumber\\
        &\quad\ \ = \lim_{r\to\infty}\V\bigg(\frac{1}{r^{\frac{d}{2}}}\sum_{n=1}^d\sum_{k\in\Z^d} \partial_n f((k+U)/r) \int x_n\, \Xi_k(\mathrm dx)\bigg)\nonumber\\
        &\quad\ \ = \lim_{r\to\infty} \frac{1}{r^d}\sum_{n,l=1}^d \sum_{s,k\in\Z^d} C_{n,l}(s,k)\int_{[0,1)^d}\partial_n f((s+u)/r)\partial_l f((k+u)/r)\, \mathrm du\nonumber\\
        &\quad\ \ = \lim_{r\to\infty} \frac{1}{r^d}\sum_{n,l=1}^d \sum_{k\in\Z^d} C_{n,l}(0,k)\int\partial_n f(s/r)\partial_l f((s+k)/r)\, \mathrm ds\nonumber\\
        &\quad\ \ = \lim_{r\to\infty}\sum_{n,l=1}^d \sum_{k\in\Z^d} C_{n,l}(0,k)\int\partial_n f(s)\partial_l f(s+k/r)\, \mathrm ds\nonumber\\
        &\quad\ \ = \sigma^2,\label{e:var-asamp}
    \end{align}
    with $\sigma^2$ defined as in \eqref{e:asymptotic-variance} and 
    \begin{equation*}
        C_{n,l}(s, k) := \C\bigg(\int x_n\, \Xi_s(\mathrm dx),\int x_l\, \Xi_k(\mathrm dx)\bigg),\quad n,l\in\{1,...,d\}, s,k\in\Z^d.
    \end{equation*}
    Therefore, $\tilde\Phi$ is $2$-uniform. Further, if $f\neq0$, then $(\langle\partial_n f, \partial_l f\rangle)_{n,l\in\{1,...,d\}}$ is a positive-definite matrix. Hence, $\sigma^2$ is the Frobenius inner product of a positive-definite matrix with a positive semi-definite matrix. Consequently, $\sigma^2=0$ is equivalent to the positive semi-definite matrix being $0$, which establishes the criterion for $\sigma^2>0$ in the proposition.

    \smallskip\noindent
    \underline{Central limit theorem:}
    
    \smallskip\noindent
    The central limit theorem is trivial for $\sigma^2=0$. For the case $\sigma^2>0$, we rely on a criterion from~\cite{JP09}. There, conditions are derived, under which
    \[
        \frac{1}{\sigma_n}\sum_{k\in D_n} Z_{k,n} \xrightarrow[\mathcal D]{n\to\infty} \mathcal N(0,1),
    \]
    where $\sigma_n := \sqrt{\V(\sum_{k\in D_n} Z_{k,n})}$ and $(Z_{k,n})_{k\in D_n, n\in\N}$ are centered random variables. Here, we first show the central limit theorem for the non-stationary model $\Phi$, and then generalize to $\tilde\Phi$. Given a sequence $(r_n)_{n\in\N}$ such that $r_n\to\infty$ as $n\to\infty$, we choose
    \begin{align*}
        D_n &:= \Z^d\cap B_{b r_n},\\
        Z_{k,n} &:= \sum_{l=1}^d \partial_l f(k/r_n)\bigg(\int x_l \,\Xi_k(\mathrm dx) - \E\bigg[\int x_l \,\Xi_k(\mathrm dx)\bigg]\bigg) ,\quad k\in D_n, n\in\N,
    \end{align*}
    where $b>0$ is chosen such that $\textnormal{supp}(f)\subseteq B_{b}$. We continue with checking the assumptions from~\cite{JP09}. Assumption 1 is satisfied since $D_n\subseteq\Z^d, n\in\N,$ and the nearest neighbor distance in $\Z^d$ is always $1$. To check Assumption 2, we choose
    \begin{equation*}
        c_{k,n} := \max_{l\in\{1,...,d\}} \|\partial_l f\|_\infty,\quad k\in D_n, n\in\N,
    \end{equation*}
    and
    \begin{equation*}
        Z := \sum_{l=1}^d\bigg|\int x_l \,\Xi_0(\mathrm dx) - \E\bigg[\int x_l \,\Xi_0(\mathrm dx)\bigg] \bigg|.
    \end{equation*}
    Then Assumption 2 is satisfied since
    \begin{equation*}
        \lim_{t\to\infty}\sup_{n\in\N}\sup_{k\in D_n} \E\big[|Z_{k,n}/c_{k,n}|^2\I\{|Z_{k,n}/c_{k,n}|>t\}\big] \leq \lim_{t\to\infty} \E\big[Z^2\I\{Z>t\}\big]=0,
    \end{equation*}
    where we used dominated convergence and that, by \eqref{e:moment-assumption},
    \begin{equation*}
        \E\bigg[\bigg(\int x_l\, \Xi_0(\mathrm dx)\bigg)^2\bigg] \leq d\E[N_0^2]<\infty, \quad l\in\{1,...,d\}.
    \end{equation*}
    Assumption 3 (b) and (c) are direct consequences of \eqref{e:mixing-condition} since the mark field $(M_k)_{k\in\Z^d}$ carries the correlations of the random field $(Z_{k,n})_{k\in D_n}$ for any $n\in\N$. The verification of Assumption 3 (a) is more involved. First, we can apply \eqref{e:mixing-condition} again to bound the function $\alpha_{\textnormal{inv}}$, defined as there, by
    \begin{equation*}
        \alpha_{\textnormal{inv}}(u) \leq c^{\frac{1}{d+\varepsilon}} u^{-\frac{1}{d+\varepsilon}}, \quad u\in(0,1).
    \end{equation*}
    Then, we apply \eqref{e:moment-assumption} to bound the upper-quantile function $\overline{Q}^{(t)}_{k,n}$, defined as there, by
    \begin{equation*}
        \overline{Q}^{(t)}_{k,n}(u) \leq \E\big[Z^m\I\{Z>t\}\big]^{\frac{1}{m}} u^{-\frac{1}{m}}<\infty,
    \end{equation*}
    for $u\in(0,1), t\in\R, k\in D_n, n\in\N$, where we used the Markov inequality. Hence, there is a bounded function $g:\R\to[0,\infty)$ such that $g(t)\to 0$ as $t\to\infty$ and
    \begin{equation*}
    \lim_{t\to\infty}\sup_{n\in\N}\sup_{k\in D_n} \int_0^1 \alpha_{\textnormal{inv}}(u)^d \overline{Q}^{(t)}_{k,n}(u)^2\, \mathrm du \leq \lim_{t\to\infty} g(t)\int_0^1 u^{-(\frac{d}{d+\varepsilon}+\frac{2}{m})}\, \mathrm du =0
    \end{equation*}
    since $\frac{d}{d+\varepsilon}+\frac{2}{m}<1$, which establishes Assumption 3 (a). Assumption 4 is redundant since Assumption 3 is fulfilled. Finally, we can verify Assumption 5 since we can show $r_n^{-d}\sigma_n^2\to\sigma^2$, as $n\to\infty$, similar to \eqref{e:var-asamp}, whereby
    \begin{equation*}
        \liminf_{n\to\infty} \frac{\sigma_n^2}{|D_n|(\max_{k\in D_n} c_{k,n})^2} = \frac{\sigma^2}{\lambda_d(B_b)(\max_{l\in\{1,...,d\}} \|\partial_l f\|_\infty)^2} >0.
    \end{equation*}

    Now, we lift this central limit theorem from $\Phi$ to $\tilde\Phi$. By \ref{e:leading-order-decomposition}, 
    \begin{equation*}
         \int f(x/r)\, \tilde\Phi(\mathrm dx) - r^d\lambda_d(f)\E[N_0] = \frac{1}{r}\sum_{n=1}^d\sum_{k\in\Z^d} \partial_n f((k+U)/r) \int x_n\, \Xi_k(\mathrm dx) + o(r^{\frac{d}{2}-1}),
    \end{equation*}
    and the right hand side (RHS) is centered since derivatives of integrable function integrate to $0$. By Slutsky's lemma it suffices to show a central limit theorem for the term on the RHS. By the previous part, we know that the term on the RHS follows a central limit theorem if we fix $U=0$. By choosing $f$ appropriately, we have also established the same central limit theorem for $U=u$ for any $u\in[0,1)^d$. Consequently, we can conclude the central limit theorem in \eqref{e:clt} from the fact that $U$ is independent of $\Phi$ since that implies $\E[g(\Phi, U)] = \int_{[0,1)^d}\E[g(\Phi, u)]\, \mathrm du$ for any bounded function $g$.
    
    \smallskip\noindent
    \underline{Generalizing the test function:}
    
    \smallskip\noindent
    In the previous parts, we have established \eqref{e:clt} for $f\in\mathcal C_c^\infty$. Here, we generalize to the setting of the assertion in two steps. First, we assume instead that $f\in L^1$ such that \eqref{e:real-space-decay} holds, and we assume that there is a $b>0$ such that $\textnormal{supp}(\hat f)\subseteq B_b$. Choose $g\in\mathcal C_c^\infty$ with $g(0)=1$ and $\lambda_d(g)\neq 0$. For $n\in\N$, define $f_n:\R^d\to\R, f_n(x):= g(x/n)(f(x)+a_n), x\in\R^d$ with $a_n$ chosen such that $\lambda_d(f_n)=\lambda_d(f)$. Then $f_n\in\mathcal C_c^\infty$. Now, let $\varphi_r$ be the characteristic function of the left hand side (LHS) of \eqref{e:clt}, and $\varphi$ that of the RHS of \eqref{e:clt}. Similarly, let $\varphi^{(n)}_{r}$ and $\varphi^{(n)}$ be the characteristic functions of the LHS and RHS of \eqref{e:clt} with respect to the test function $f_n$ instead of $f$. We obtain that
    \begin{equation}\label{e:cf-inequality}
        |\varphi_r(k)-\varphi(k)| \leq |\varphi_r(k)-\varphi^{(n)}_r(k)| + |\varphi^{(n)}_r(k)-\varphi^{(n)}(k)| + |\varphi^{(n)}(k)-\varphi(k)|
    \end{equation}
    for $k\in\R, r>0, n\in\N$. We can apply Proposition~\ref{p:clt} with the test function $f_n\in\mathcal C_c^\infty$ instead of $f$ to see that the middle term converges to $0$ as $r\to\infty$ for any fixed $n\in\N$. To treat the left term consider
    \begin{equation*}
        (\sigma_r^{(n)})^2:=\V\bigg(\frac{1}{r^{\frac{d}{2}-1}}\int f(x/r) - f_n(x/r)\, \tilde\Phi(\mathrm dx)\bigg) = \frac{r^{d+2}}{(2\pi)^d}\int |\hat{f}(rk)-\hat{f}_n(rk)|^2\, \hat\beta_{\tilde\Phi}(\mathrm dk).
    \end{equation*}
    Since $\hat f$ is differentiable by \eqref{e:real-space-decay} and satisfies $\textnormal{supp}(\hat f')\subseteq B_b$, we can calculate that
    \begin{align*}
        |\hat{f}(rk)-\hat{f}_n(rk)| &= \frac{1}{(2\pi)^d}\bigg|\int \hat{g}(s) \big(\hat{f}(rk) -\hat{f}(rk-s/n)\big) \,\mathrm ds + a_nn^d \hat g(rnk)\bigg|\\
        &\leq \frac{\|\hat f '\|_\infty}{(2\pi)^dn}\int |\hat g(s)| \|s\|\,\mathrm ds + a_nn^d\|\hat g\|_\infty,
    \end{align*}
    and for $k\in B_{\frac{b}{r}}^c$, we obtain the alternate bound
    \begin{align*}
        |\hat{f}(rk)-\hat{f}_n(rk)| &\leq \frac{\|\hat f\|_\infty}{(2\pi)^d}  \int_{B_{nb}+nrk} |\hat g(s)|\,\mathrm ds + |a_n|n^d |\hat g(rnk)|\\
        &\leq  \frac{\|\hat f\|_\infty}{(2\pi)^dn^{d+2}(r\|k\|-b)^{d+2}}\int |\hat g(s)|\|s\|^{d+2}\,\mathrm ds + |a_n|n^d \frac{|\hat g(rnk)|(rn\|k\|)^{d+2}}{n^{d+2}(r\|k\|)^{d+2}}.
    \end{align*}
    Since $g\in\mathcal C_c^\infty$, $\hat g$ decays quicker than any polynomial and since $\tilde\Phi$ is $2$-uniform, there are $c_1,c_2>0$ independent of $r$ and $n$ such that
    \begin{align*}
        (\sigma_r^{(n)})^2 &\leq \big(\frac{c_1}{n^2}+a_n^2n^{2d}\big)\bigg( \frac{\hat\beta_{\tilde\Phi}\big(B_{\frac{b}{r}}\big)}{(\frac{b}{r})^{d+2}} + \sum_{l=2}^\infty \frac{\hat\beta_{\tilde\Phi}\big(B_{\frac{lb}{r}}\big)}{(\frac{lb}{r})^{d+2}} \frac{l^{d+2}}{(l-1)^{2d+4}}\bigg)\\
        &\leq c_2(n^{-2} + a_n^2n^{2d}).
    \end{align*}
    Further, by the choice of $a_n$, we have $\E[\int f(x/r) - f_n(x/r)\, \tilde\Phi(\mathrm dx)]=0$, so we can apply the Markov inequality to see
    \begin{equation*}
        \frac{|\varphi_r(k)-\varphi^{(n)}_r(k)|}{\|k\|} \leq \E\bigg[\bigg|\frac{1}{r^{\frac{d}{2}-1}}\int f(x/r) - f_n(x/r)\, \tilde\Phi(\mathrm dx)\bigg|\bigg] \leq \sigma_r^{(n)} \leq \sqrt{c_2(n^{-2} + a_n^2n^{2d})}.
    \end{equation*}
    In total, we can now conclude that
    \begin{equation}\label{e:cf-inequality2}
        \limsup_{r\to\infty}|\varphi_r(k)-\varphi(k)| \leq \|k\|\sqrt{c_2(n^{-2} + a_n^2n^{2d})} + |\varphi^{(n)}(k)-\varphi(k)|.
    \end{equation}
    A simple application of Fubinis theorem yields that $a_nn^d\to 0$ as $n\to\infty$, whereby it only remains to show that the right term of the bound goes to $0$ as $n\to\infty$ to prove the central limit theorem in this case. Since $\varphi^{(n)}$ and $\varphi$ are both characteristic functions of centered Gaussian random variables, it suffices to show that the variance converges. By the specific formula \eqref{e:asymptotic-variance} for the variances, the problem can be further reduced to showing $f_n'\to f'$ in $L^2$ as $n\to\infty$. This limit holds since for any $l\in\{1,...,d\}$,
    \begin{equation*}
        \|\partial_l f_n-\partial_l f\|_2 \leq \frac{1}{n} \|\partial_l g(\cdot/n) f\|_2 + |a_n|n^{\frac{d}{2}-1} \|\partial_l g\|_2 + \| (g(\cdot/n)-1) \partial_l f\|_2\xrightarrow{n\to\infty} 0,
    \end{equation*}
    where we use that $f,f'\in L^2$ and $g\in\mathcal C_c^\infty$ with $g(0)=1$. This concludes the first part of the proof. 

    In the second step, we generalize $f$ to the setting of the assertion by assuming \eqref{e:fourier-space-decay} instead of the existence of a $b>0$ such that $\textnormal{supp}(\hat f)\subseteq B_b$. Without loss of generality, we assume $\vartheta\in(0,1)$. Also, we replace $g$ by a Schwartz function that satisfies $\lambda_d(g)=1$ and $\hat g\in\mathcal C_c^\infty$, and we instead define $f_n(x):= \int g(y) f(x-y/n)\, dy, x\in\R^d$. Then $f_n$ has all properties that were assumed for $f$ in the first step of this part. When we redefine $\varphi_r, \varphi^{(n)}_r, \varphi^{(n)}, \varphi$ with respect to the new functions $f$ and $f_n$ and consider \eqref{e:cf-inequality}, the middle term vanishes as $r\to\infty$ by the first part of the proof. To bound the left term, we bound
    \begin{equation*}
        |\hat{f}(rk)-\hat{f}_n(rk)| = |1-\hat g(rk/n)| |\hat{f}(rk)| \leq c_3 \min((r\|k\|/n)^{\vartheta/4}, 1) |\hat{f}(rk)|,
    \end{equation*}
    where $c_3>0$ exists since $g\in\mathcal S(\mathbb R^d)$. Since $\hat f$ decays sufficiently quickly, we derive similarly to before that, with $\sigma_r^{(n)}$ defined in terms of the new $f$ and $f_n$, 
    \begin{equation*}
        \limsup_{r\to\infty}\sigma_r^{(n)} \leq \frac{c_4}{n^{\vartheta/4}},
    \end{equation*}
    where $c_4>0$ is independent of $n$. We can then derive an analogue of \eqref{e:cf-inequality2} and simplify the remainder of the assertion to showing $f_n'\to f'$ in $L^2$ as $n\to\infty$. Finally, by Plancherels theorem and dominated convergence,
    \begin{align*}
        \sum_{l=1}^d\|\partial_l f_n-\partial_l f\|_2^2 &= \frac{1}{(2\pi)^{d}} \int \|k\|^2|\hat{f}(k)-\hat{f}_n(k)|^2 \,\mathrm dk\\
        &= \frac{1}{(2\pi)^{d}} \int |1-\hat g(k/n)|^2 \|k\|^2 |\hat{f}(k)|^2 \,\mathrm dk \\
        &\xrightarrow{n\to\infty} 0,
    \end{align*}
    since $\hat{g}(0)=1$ and $\hat{f}$ decays sufficiently quickly. This fact concludes the proof.

\end{proo}

\section{Examples of decorated lattices}\label{sec:6}

The following examples illustrate how cell counts, dependence between cells, and local point arrangements contribute to
the spectrum of a decorated lattice.
We first recover familiar models, including Poisson processes and independently perturbed lattices. We then compare independent and dependent cell counts, showing how the random environment can produce hyperuniformity
even when the number of points varies between cells. Finally, we examine how random within-cell distributions and constraints on motif moments affect the finer order of uniformity and the associated fluctuations.

\subsection{Poisson motifs}

Let $\mu$ be a finite measure on $Q$, and let the motifs
$(\Xi_j)_{j\in\mathbb Z^d}$ be independent Poisson processes
with common intensity measure $\mu$.
This is a decorated lattice with a singleton mark space. Conditional on $U$, the stationarized process is Poisson
with intensity measure
\[
\Lambda_U
=
\sum_{j\in\mathbb Z^d}\theta_{j+U}\mu.
\]
Thus, in general, $\widetilde\Phi$ is a Cox process directed
by a randomly shifted periodic measure.

Theorem~\ref{L_Bartlett_CIM} gives
\[
\widehat\beta(\mathrm dk)
=
\mu(Q)\,\mathrm dk
+
(2\pi)^d \sum_{z\in\mathbb Z^d\setminus\{0\}}
|\widehat\mu(2\pi z)|^2
\delta_{2\pi z}(\mathrm dk).
\]
In particular, if $\mu=\gamma\lambda_d|_Q $,
then
$\widehat\mu(2\pi z)=0 $ for every $z\neq0.$
The resulting process is a homogeneous Poisson process of intensity
$\gamma$, and its Bartlett spectral measure is
\[
\widehat\beta(\mathrm dk)
=
\gamma\,\mathrm dk.
\]
which is in correspondence with \eqref{Eq_Bartlett_Poisson}.
If $\mu(Q)>0$, then $S(0)=\mu(Q)>0$, so the process is of course
not hyperuniform.

\subsection{Independent motifs}

When the mark field is i.i.d., the unconditional motifs
are also i.i.d. Let $\Psi$ have their common law and
assume $\E[\Psi(Q)^2]<\infty$.
The continuous part of the Bartlett spectrum has density
\[
S(k)
=
\V\left(\int_Q e^{-ik\cdot x}\,\Psi(\mathrm dx)\right).
\]
In particular, $S(0)=\V(\Psi(Q))$, so hyperuniformity
requires and is implied by a deterministic motif count.
The following example recovers the familiar one-point case.
\begin{exam}[Independently perturbed lattice with bounded perturbations] 
Let $(Y_k)_{k\in\mathbb Z^d}$ be i.i.d.\ random variables taking
values in $Q$, with common distribution $\vartheta$, and define
\[
\Phi
:=
\sum_{k\in\mathbb Z^d}\delta_{k+Y_k}
\]
and $\widetilde\Phi$ its stationarized version.
This model is a special case of the decorated lattice with the
singleton mark space $\mathbb M=\{m\}$
and motif law 
$P_m=\mathcal L(\delta_Y),
 Y\sim\vartheta$.
In particular,
\[
n(m)=s(m)=1,
\]
and the mark-covariance measures satisfy $\eta_r=0, r\in\mathbb Z^d.$
Hence
$\Gamma_r(k)=0$ for every $r\in\mathbb Z^d.$

The motif intensity measure is
$\mu_m=\vartheta,$
whereas its second-order moment measure is concentrated on the
diagonal. 
Consequently,
$$\widehat M_m^{(2)}(k)=1 \qquad \text{and} \qquad \widehat\mu_m(k)
=
\int_Qe^{-ik\cdot x}\,\vartheta(\mathrm dx)=: \varphi_Y(k)$$
is the characterictic function of $Y$.
Theorem~\ref{L_Bartlett_CIM} therefore yields
$$
\widehat\beta_{\widetilde\Phi}(\mathrm dk)
=
\left(
1-\left|\varphi_Y(k)\right|^2
\right)\mathrm dk
+
(2\pi)^d \sum_{z\in\mathbb Z^d\setminus\{0\}}
\left|\varphi_Y(2\pi z)\right|^2
\delta_{2\pi z}(\mathrm dk),
$$
which is in correspondence with Lemma 2 of \cite{F26}. Further, Proposition~\ref{p:clt} yields that $\tilde\Phi$ is $2$-uniform and satisfies the central limit theorem
\begin{equation*}
    \frac{1}{r^{\frac{d}{2}-1}} \bigg(\int (2\pi)^{-\frac{d}{2}}e^{-\frac{x^{\scalebox{0.7}{$\scriptstyle T$}}\cdot x}{2r^2}}\, \tilde\Phi(\mathrm dx)-r^d\bigg) \xrightarrow{\mathcal D} \mathcal N(0, \sigma^2),
\end{equation*}
with
\begin{equation*}
    \sigma^2 = \frac{\textnormal{tr}(\C(Y))}{2^{d+1}\pi^\frac{d}{2}}.
\end{equation*}
These results coincide with those obtained in~\cite{Ma25}. 
In particular, if $\vartheta$ is the uniform distribution on $Q$, then
\[
\varphi_Y(2\pi z)=0
\qquad
\text{for every }z\in\mathbb Z^d\setminus\{0\},
\]
and hence the Bartlett spectral measure is purely absolutely
continuous:
\[
\widehat\beta_{\widetilde\Phi}(\mathrm dk)=
\left[
1-
\prod_{j=1}^d
\left(
\frac{\sin(k_j/2)}{k_j/2}
\right)^2
\right]\mathrm dk,
\]
where the expression is understood as continuous when $k_j=0$. The disappearance of the nonzero Bragg peaks is known
as cloaking (see \cite{KKT20} or \cite[Section 10.3]{C21}).
It removes the discrete part of the spectrum while preserving hyperuniformity. In this particular case, the asymptotic variance in the central limit theorem from above is
\begin{equation*}
    \sigma^2 = \frac{d}{2^{d+3}3\pi^\frac{d}{2}}.
\end{equation*}
\end{exam}

\begin{figure}[t]
    \centering
    \includegraphics[width=\textwidth]{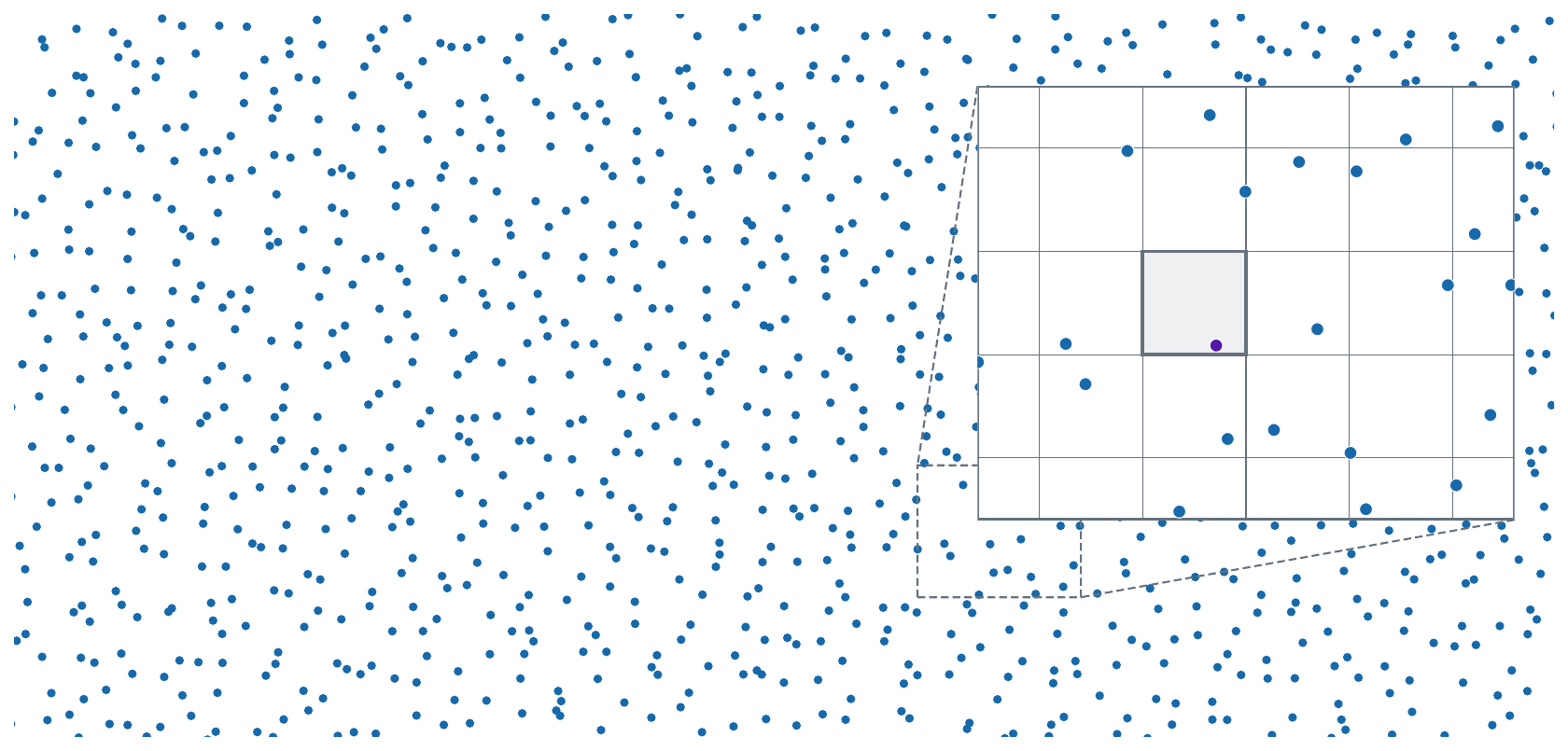}
    \caption{Sample of a stationarized decorated lattice with one point per cell. Here, the single point is distributed uniformly in its cell and its relative location is independent of the relative position of points in other cells. Hence, the stationarized decorated lattice is an i.i.d. perturbed lattice with uniformly bounded perturbations. Therefore, it is $2$-uniform and not beyond $2$-uniform. In this particular case, the uniform distribution yields that the Bartlett spectral measure of the stationarized decorated lattice is absolutely continuous, whereby it is also called a cloaked lattice, see~\cite{KKT20}.}
    \label{fig:cloaked_lattice}
\end{figure}

\subsection{Binomial motifs}
In this example, we consider binomial motifs and investigate three settings: a fixed count in every
cell, independent two-valued counts, and two-valued counts balanced
across neighbouring pairs. In all three cases, points within a cell
are sampled independently from the same law $\vartheta$, given the
cell count. The last two examples show that even with the same
one-cell count distribution, independent fluctuations destroy
hyperuniformity, while balanced pairs restore it.

Let $\vartheta$ be a probability measure on $Q$.
For every $n\in\mathbb N_0$, let $P_n$ be the law of the binomial point process 
\[ \Xi^{(n)} = \sum_{i=1}^{n}\delta_{Y_i}, \] where $Y_1,\ldots,Y_n$ are independent with common distribution $\vartheta$. Thus a mark prescribes exactly $n$ points in its cell. Let $(M_j)_{j\in\mathbb Z^d}$ be a stationary $\mathbb N_0$-valued mark field. Conditionally on the mark field, let the motifs $(\Xi_j)_{j\in\mathbb Z^d}$ be independent and satisfy 
\[ \mathcal L(\Xi_j\mid\mathbf M)=P_{M_j}. \]

For a mark $n\in\mathbb N_0$, the mean motif measure is \( \mu_n=n\vartheta, \) and hence 
\(\widehat\mu_n(k) = n\varphi(k) \), where $\varphi$ is the characteristic function of $\vartheta$.
Moreover, 
\[\widehat M_n^{(2)}(k) = \mathbb E\left[ \left| \sum_{i=1}^{n}e^{-ik\cdot Y_i} \right|^2 \right] = n+n(n-1)|\varphi(k)|^2. \]
Consequently, 
\[
 \widehat M_n^{(2)}(k) - |\widehat\mu_n(k)|^2 = n\left(1-|\varphi(k)|^2\right). 
 \] 
Put $ \overline m:=\mathbb E[M_0]. $ For the mark-covariance term \eqref{EQ_Gamma_CIM}, one has 
\[
\Gamma_r(k) = \C\left( \widehat\mu_{M_0}(k), \widehat\mu_{M_r}(k) \right) = |\varphi(k)|^2 \C(M_0,M_r). 
\]
Therefore, whenever the covariance series is absolutely summable and $\E[M_0^2]<\infty$, 
\begin{equation} S(k) = \overline m\left(1-|\varphi(k)|^2\right) + |\varphi(k)|^2 \sum_{r\in\mathbb Z^d} e^{ik\cdot r} \C(M_0,M_r). \label{EQ_binomial_general_spectral_density} \end{equation} 
Furthermore, 
\( \mathbb E\widehat\mu_M(k) = \overline m\,\varphi(k), \)
and hence Theorem~\ref{L_Bartlett_CIM} gives 
\begin{align*} \widehat\beta(\mathrm dk) ={}& S(k)\,\mathrm dk \notag+ (2\pi)^d\overline m^{\,2} \sum_{z\in\mathbb Z^d\setminus\{0\}} |\varphi(2\pi z)|^2 \delta_{2\pi z}(\mathrm dk). \end{align*} 

\begin{remark}[Uniformly distributed points within cells]
Suppose that $\vartheta=\lambda_d|_Q$. Then
\[
|\varphi(k)|^2
=
\prod_{j=1}^d
\left(\frac{\sin(k_j/2)}{k_j/2}\right)^2
,
\]
with each factor interpreted as one at $k_j=0$.
In particular, $\varphi(2\pi z)=0$ for every
$z\in\mathbb Z^d\setminus\{0\}$, so all nonzero Bragg
peaks disappear.
Thus, uniform sampling removes the discrete part of the
spectrum.
\end{remark}

\subsection*{Common number of points in every cell} 
Fix $n\in\mathbb N$. Let 
\( M_j=n\) for every $j\in\mathbb Z^d.$
This may be regarded as a degenerate i.i.d.\ mark field. Every cell contains exactly $n$ independent points with common distribution $\vartheta$. Since the mark field is deterministic, it is also uncorrelated and
the Bartlett spectral measure is 
\[
\widehat\beta(\mathrm dk) = n\left(1-|\varphi(k)|^2\right)\mathrm dk + (2\pi)^d n^2 \sum_{z\in\mathbb Z^d\setminus\{0\}} |\varphi(2\pi z)|^2 \delta_{2\pi z}(\mathrm dk). 
\]
In particular,
\(
S(k) = n\left(1-|\varphi(k)|^2\right). 
\) Since $S(0)=0$ and the process is hyperuniform. 

Let  $\Sigma:=\C(Y).$
The random vector $Y$ is bounded and therefore, 
\[ 1-|\varphi(k)|^2 = k^{\mathsf T}\Sigma k+o(|k|^2), \qquad k\to0. \]
Writing 
\( \kappa_d:=\lambda_d(B_1), \)
we obtain 
\[
\frac{ \widehat\beta(B_\varepsilon) }{ \varepsilon^d } = \frac{n\kappa_d}{d+2} \operatorname{tr}(\Sigma)\, \varepsilon^2 + o(\varepsilon^2). \]
Hence, if $\vartheta$ is not concentrated at a single point, the process is $2$-uniform but not beyond $2$-uniform. Proposition~\ref{p:clt} yields the central limit theorem
\begin{equation*}
    \frac{1}{r^{\frac{d}{2}-1}} \bigg(\int (2\pi)^{-\frac{d}{2}}e^{-\frac{x^{\scalebox{0.7}{$\scriptstyle T$}}\cdot x}{2r^2}}\, \tilde\Phi(\mathrm dx)-nr^d\bigg) \xrightarrow{\mathcal D} \mathcal N(0, \sigma^2),
\end{equation*}
with
\begin{equation*}
    \sigma^2 = \frac{\textnormal{tr}(\Sigma)n}{2^{d+1}\pi^\frac{d}{2}}.
\end{equation*}

We now keep the number of points fixed and let the environment select their distribution within each cell. This separates the effect of the mark field on local geometry from its effect on cell counts. Let $P_\vartheta$ be the law of the binomial point process
\begin{equation*}
    \Xi^{(\vartheta)} = \sum_{l=1}^n\delta_{Y_l},
\end{equation*}
where $Y_1,...,Y_n$ are independent with common distribution $\vartheta$. Let $(M_j)_{j\in\Z^d}$ be a stationary mark field of probability measures on $Q$ such that
\[
\sum_{r\in\mathbb Z^d}|\eta_r|(\mathbb M^2)<\infty.
\]
Since every motif contains exactly $n$ points, this is
the summability assumption of
Theorem~\ref{L_Bartlett_CIM} up to the constant factor $n^2$. Then
\[
S(k)
=
n\left(1-\E[|\widehat M_0(k)|^2]\right)
+
n^2\sum_{r\in\mathbb Z^d}
e^{ik\cdot r}
\C\bigl(\widehat M_0(k),\widehat M_r(k)\bigr),
\]
and
\[
\widehat\beta_{\widetilde\Phi}(\mathrm dk)
=
S(k)\,\mathrm dk
+
(2\pi)^d n^2
\sum_{z\in\mathbb Z^d\setminus\{0\}}
\left|\E[\widehat M_0(2\pi z)]\right|^2
\delta_{2\pi z}(\mathrm dk).
\]
Since $\widehat M_j(0)=1$ almost surely, $S(0)=0$.
Thus, under these assumptions, the process is hyperuniform.
The mark field changes the local distributions and their
dependence, while preserving the count in every cell.

Under the additional assumption that the mark field satisfies the $\alpha$-mixing condition \eqref{e:mixing-condition} for some $\varepsilon>0$, we can conclude from Proposition~\ref{p:clt} that the decorated lattice $\tilde\Phi$ is\linebreak$2$-uniform and satisfies the central limit theorem
\begin{equation*}
    \frac{1}{r^{\frac{d}{2}-1}} \bigg(\int (2\pi)^{-\frac{d}{2}}e^{-\frac{x^{\scalebox{0.7}{$\scriptstyle T$}}\cdot x}{2r^2}}\, \tilde\Phi(\mathrm dx)-nr^d\bigg) \xrightarrow{\mathcal D} \mathcal N(0, \sigma^2),
\end{equation*}
where
\begin{equation*}
    \sigma^2 = \frac{1}{2^{d+1}\pi^\frac{d}{2}} \Bigg(n \E\big[\textnormal{tr}(V_{M_0})\big]+n^2\sum_{k\in\Z^d} \textnormal{tr}\bigg(\C\bigg(\int x\, M_0(\mathrm dx), \int x\, M_k(\mathrm dx)\bigg)\bigg)\Bigg),
\end{equation*}
with $V_\vartheta$ being the variance of a random vector with distribution $\vartheta$. Note that both summands are nonnegative. Here, one can observe that the scaling of the asymptotic variance in terms of $n$ changes when the random field of the centroids of the marks is hyperuniform and the right term drops out. This is connected to the fact that the asymptotic variance of $\frac{1}{n}\tilde\Phi$ approaches that of a stationarized lattice which is perturbed by the centroids of the mark field. A perturbed lattice is not only $2$-uniform but also beyond $2$-uniform if the displacement field is hyperuniform.


\subsection*{Independent two-valued cell counts}
Fix two distinct integers $ m_1,m_2\in\mathbb N_0 $ and let $(M_j)_{j\in\mathbb Z^d}$ be i.i.d.\ with 
\( \mathbb P(M_j=m_1)=p,\mathbb P(M_j=m_2)=1-p, \)
where $p\in(0,1)$. Put 
\[ \overline m = pm_1+(1-p)m_2 \] 
and 
\[ \sigma_M^2 := \V(M_0) = p(1-p)(m_1-m_2)^2. \]
Formula~\eqref{EQ_binomial_general_spectral_density} therefore gives 
\[
S(k) = \overline m\left(1-|\varphi(k)|^2\right) + \sigma_M^2|\varphi(k)|^2. 
\]
At the origin, 
$ S(0) = \sigma_M^2 >0. 
$
Thus the process is not hyperuniform. More precisely, 
\[ \frac{ \widehat\beta(B_\varepsilon) }{ \varepsilon^d } \longrightarrow \kappa_d\sigma_M^2, \qquad \varepsilon\downarrow0. \] 
Hence the process is $0$-uniform but not beyond $0$-uniform. The intrinsic motif term vanishes at the origin because the number of points is fixed after the mark is known. The failure of hyperuniformity is caused entirely by the independent fluctuations of the mark-selected cell counts. Since $\tilde\Phi$ has these Poisson-like fluctuations, we also obtain a classical central limit theorem from Proposition \ref{p:classical-clt}, i.e., 
\begin{equation*}
    r^{-\frac{d}{2}}\big(\tilde\Phi(B_r) - \overline{m}\kappa_d r^d\big) \xrightarrow{\mathcal D} \mathcal N(0, \kappa_d\sigma_M^2).
\end{equation*}

Independence of the cell counts is decisive in the preceding example. We now retain the same two possible counts but introduce dependence that makes neighboring fluctuations compensate.

\subsection*{Dependent two-valued counts with balanced pairs} Fix two distinct integers $ m_1,m_2\in\mathbb N_0. $
We retain the one-site count distribution of the preceding
example with $p=1/2$, but replace independence by a balanced-pair construction. Let 
\[ e_1=(1,0,\ldots,0)\in\mathbb Z^d \]
and let $R$ be uniformly distributed on $\{0,1\}$. Conditionally on $R$, partition $\mathbb Z^d$ into the pairs 
\[ \left\{ (2\ell+R,v), (2\ell+R+1,v) \right\}, \qquad \ell\in\mathbb Z,\quad v\in\mathbb Z^{d-1}. \]
Independently for different pairs, assign the ordered marks either as 
\[ (m_1,m_2) \qquad\text{or}\qquad (m_2,m_1), \]
each with probability $1/2$. The random phase $R$ makes the resulting mark field stationary under translations of $\mathbb Z^d$. Every site satisfies 
\[ \mathbb P(M_j=m_1) = \mathbb P(M_j=m_2) = \frac12, \]
and therefore 
\( \overline m = \frac12 (m_1+m_2). \) 
Put $\Delta:=m_1-m_2.$ The autocovariance function of the mark field is 
\[
 \C(M_0,M_r) = \begin{cases} \displaystyle \frac{\Delta^2}{4}, &r=0,\\[0.6em] \displaystyle -\frac{\Delta^2}{8}, &r=\pm e_1,\\[0.6em] 0, &\text{otherwise}. \end{cases} 
\]
In particular, 
\[
\sum_{r\in\mathbb Z^d} \C(M_0,M_r) =  0. 
\]
The lattice spectral density of the cell-count field is 
\[
\sum_{r\in\mathbb Z^d} e^{ik\cdot r} \C(M_0,M_r) = \frac{\Delta^2}{4} -\frac{\Delta^2}{8}e^{ik_1} -\frac{\Delta^2}{8}e^{-ik_1} = \frac{\Delta^2}{4}(1-\cos k_1). \]
Consequently, 
\[ S(k) = \overline m\left(1-|\varphi(k)|^2\right) + \frac{\Delta^2}{4} |\varphi(k)|^2(1-\cos k_1) \]
and in particular,  $S(0)=0.$ 
Thus, although $\V(M_0)=\Delta^2/4>0$, the process is
hyperuniform. Each pair has the deterministic total count
$m_1+m_2$, so fluctuations in neighboring cells compensate
exactly. This illustrates how the random environment can
produce hyperuniformity without fixing the count in
each individual cell.

\begin{figure}[h]
    \centering
    \includegraphics[width=\textwidth]{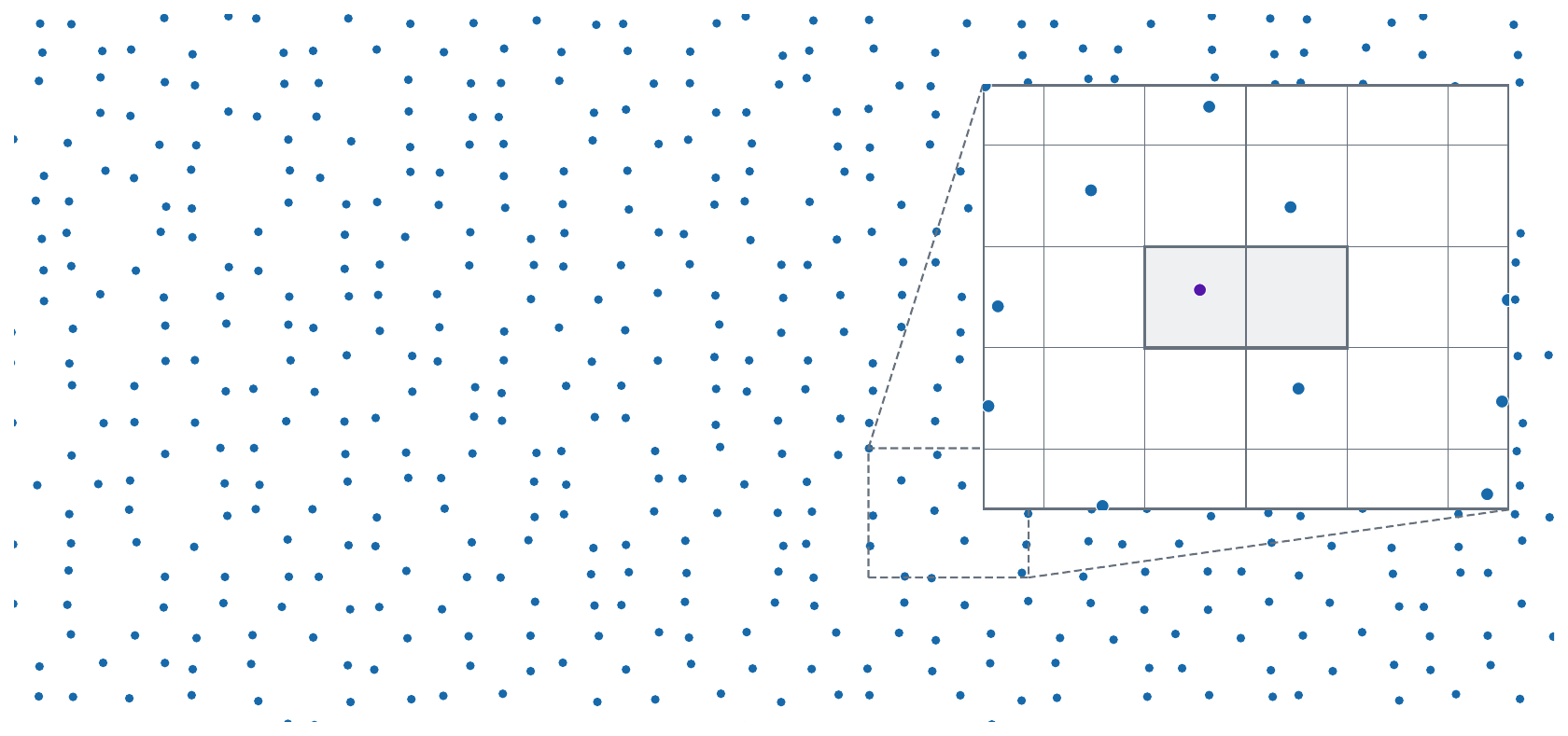}
    \caption{Sample of a stationarized decorated lattice with binomial motifs with dependent two-valued counts and balanced pairs. The underlying lattice can be divided into pairs of neighboring cells, such that one of them contains $m_1=1$ and the other contains $m_2=0$ points. Both permutations have equal probability and are chosen independently for different pairs. If there is a point in the cell, it is the centroid that is displaced by a random vector following $\mathcal U([-0.15, 0.15)^2)$ independent of the mark field and the displacements in other cells. Hence, the stationarized decorated lattice is a thinning of an i.i.d. perturbed lattice that is $2$-uniform and not beyond $2$-uniform.}
    \label{fig:binomail_balanced_pairs}
\end{figure}

Let $\Sigma:=\C(Y)$ for $Y\sim\vartheta$.
As $k\to0$,
\[
S(k)
=
\overline m\,k^{\mathsf T}\Sigma k
+
\frac{\Delta^2}{8}k_1^2
+
o(|k|^2).
\]
Consequently,
\[
\frac{\widehat\beta_{\widetilde\Phi}(B_\varepsilon)}
{\varepsilon^d}
=
\frac{\kappa_d}{d+2}
\left(
\overline m\,\operatorname{tr}(\Sigma)
+\frac{\Delta^2}{8}
\right)\varepsilon^2
+
o(\varepsilon^2).
\]
Since $\Delta\neq0$, the process is $2$-uniform but not
beyond $2$-uniform, even when the point locations within
each cell are deterministic. Proposition~\ref{p:clt} does not apply directly.
The count field is not beyond $2$-uniform, and the common
pairing phase prevents the required mixing condition.
Nevertheless, the construction becomes an independent
motif model when each pair is regarded as a single cell
of volume two. This provides a separate route to a CLT.
Grouping each balanced pair into a cell
$D=[0,2)\times[0,1)^{d-1}$ gives i.i.d.\ motifs with
deterministic count $m_1+m_2$. Their first-moment vector
$W$ satisfies
\[
\C(W)
=
2\overline m\,\Sigma
+\frac{\Delta^2}{4}e_1e_1^{\mathsf T},
\qquad \Sigma:=\C(Y),\quad Y\sim\vartheta.
\]
A first-order Taylor expansion and the Feller-Lindeberg central limit theorem therefore give, for real-valued
$f\in\mathcal S(\mathbb R^d)$,
\[
r^{1-d/2}
\left(
\int f(x/r)\,\widetilde\Phi(\mathrm dx)
-\overline m\,r^d\int f
\right)
\xrightarrow{\mathcal D}
\mathcal N(0,\sigma_f^2),
\]
with
\[
\sigma_f^2
=
\overline m\int\nabla f(x)^{\mathsf T}\Sigma\nabla f(x)\,\mathrm dx
+\frac{\Delta^2}{8}\|\partial_1f\|_2^2.
\]
The argument is applied conditional on the uniform shift
$U+Re_1$. Averaging preserves the same Gaussian limit.
Schwartz decay makes the shift-dependent mean correction
negligible. For $f(x)=(2\pi)^{-d/2}e^{-|x|^2/2}$, this reduces to
\[
\sigma_f^2
=
\frac{\overline m\,\operatorname{tr}(\Sigma)+\Delta^2/8}
{2^{d+1}\pi^{d/2}}.
\]

\subsection{Deterministic motifs and higher-order uniformity}

In this example we show how constraints on the spatial moments of
independently selected motifs produce higher orders
of uniformity.
Each mark selects a deterministic configuration, and
the first moment order that varies between configurations
determines the leading spectral order.

Let $(\mathbb M,\mathcal M)$ be a standard Borel
space and let
\[
m\longmapsto \nu_m\in\mathcal N_f(Q)
\]
be a measurable family of deterministic finite point measures on
$Q=[0,1)^d$. Let $(M_j)_{j\in\mathbb Z^d}$ be an i.i.d.\ mark field
with common law $\pi$, and set
\[
\Xi_j=\nu_{M_j}.
\]

In the context of \cite{LK26}, we can view the stationarized decorated lattice $\tilde\Phi$ as the destination of a random transport that originates in the stationarized lattice $\Z^d+U$. In particular,
\begin{equation*}
    \tilde\Phi = \int \nu_{M_{x-U}}(\cdot -x)\, (\Z^d+U)(\mathrm dx).
\end{equation*}
Since the family $(\nu_{M_j})_{j\in\Z^d}$ inherits independence from the mark field, this transport satisfies the assumptions of \cite[Theorem 5.12]{LK26}. Hence, the following lemma is a special case of that theorem. Given that $\nu_{M_0}(Q)$ is deterministic, it is an even more direct special case of \cite[Example 5.15]{LK26}. For the convenience of the reader, we still give a direct proof here. Another construction of hyperuniform cluster processes from lattices, viewed from a slightly different perspective, can be found in \cite[Subsection 2.5]{La25}.

The level of uniformity can be classified based on the following quantities. For a multi-index $\alpha\in\mathbb N_0^d$, define
\[
m_\alpha(m):=\int_Q x^\alpha\,\nu_m(\mathrm dx).
\]

\begin{lemm}\label{l:wang}
Assume that
\begin{equation}\label{EQ_second_number_counts}
    \int_{\mathbb M}\nu_m(Q)^2\,\pi(\mathrm dm)<\infty.
\end{equation}
Suppose that, for some $q\geq1$, the moments of all orders below $q$
are independent of the mark, namely there are constants $c_\alpha, |\alpha|<q$ such that
\[
m_\alpha(m)=c_\alpha
\qquad
\text{for $\pi$-a.e. }m\in\mathbb M
\text{ and all }|\alpha|<q.
\]
Assume moreover that the order-$q$ moments are not all
$\pi$-a.s.\ constant. 
Then, the stationarized decorated lattice is
$2q$-uniform but not beyond $2q$-uniform.
\end{lemm}

\begin{figure}[t]
    \centering
    \includegraphics[width=\textwidth]{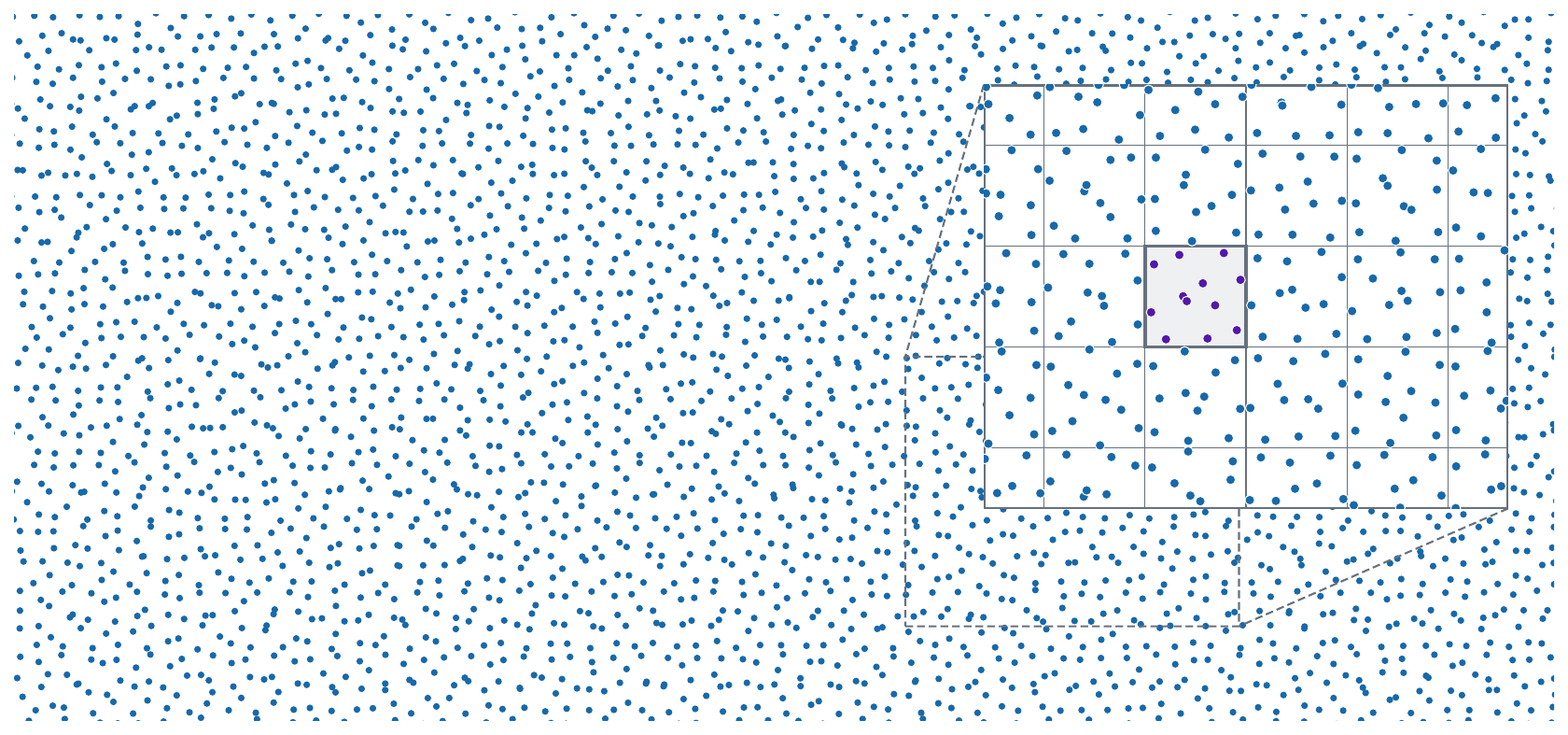}
    \caption{Sample of a stationarized decorated lattice with $c_{(0,0)}=12$ points per cell chosen such that all normalized moments of order less than $q=6$ are equal to the moments of the cell. The points in each cell are obtained from collective coordinate optimization utilizing Newton's algorithm with random, independent initialization as in \cite{LK26}. By Lemma \ref{l:wang}, the stationarized decorated lattice is $12$-uniform. We have strong numerical evidence that this procedure yields random $6$th-order moments, whereby the stationarized decorated lattice is not beyond $12$-uniform by Lemma \ref{l:wang}.}
    \label{fig:moment_matching}
\end{figure}

The order of uniformity is determined by the first spatial
moment that varies between motifs. Equal point counts remove
zeroth-order fluctuations; equal first moments also remove
the quadratic spectral contribution. More generally, fixing
all moments below order $q$ postpones the first nonzero
spectral term to order $2q$.

\begin{proo}
Since the motifs are deterministic once the mark is specified, there is
no intrinsic motif randomness. Moreover, since the mark field is
i.i.d., there are no covariance contributions between distinct cells.
Hence the absolutely continuous spectral density is
\[
S(k)
=
\operatorname{Var}\bigl(\widehat\nu_M(k)\bigr).
\]

By Taylor expansion,
\[
\widehat\nu_m(k)
=
\sum_{|\alpha|<q}
\frac{(-i)^{|\alpha|}}{\alpha!}
m_\alpha(m)k^\alpha
+
(-i)^q
\sum_{|\alpha|=q}
\frac{m_\alpha(m)}{\alpha!}k^\alpha
+
R_m(k),
\]
where, since $Q$ is bounded,
\[
|R_m(k)|
\leq C_{Q,q}\,\nu_m(Q)|k|^{q+1}.
\]
Under the assumption \eqref{EQ_second_number_counts},
we have
\[
\left(\E|R_M(k)|^2\right)^{1/2}\leq C |k|^{q+1}.
\]
The moments of order $<q$ do not depend on the mark. Hence they disappear after centering. If
\[
\overline m_\alpha
:=
\int_{\mathbb M}m_\alpha(m)\,\pi(\mathrm dm),
\qquad |\alpha|=q,
\]
then
\[
\widehat\nu_M(k)-\mathbb E[\widehat\nu_M(k)]
=
(-i)^q
\sum_{|\alpha|=q}
\frac{m_\alpha(M)-\overline m_\alpha}{\alpha!}
k^\alpha
+
O_{L^2}(|k|^{q+1}).
\]
Taking the $L^2$-norm squared gives
\[
S(k)
=
H_{2q}(k)+o(|k|^{2q}),
\]
where
\[
H_{2q}(k)
:=
\mathbb E\left[
\left|
\sum_{|\alpha|=q}
\frac{m_\alpha(M)-\overline m_\alpha}{\alpha!}
k^\alpha
\right|^2
\right].
\]
The function $H_{2q}$ is nonnegative and 
\(H_{2q}(t k) = t^{2q}H_{2q}(k), t>0.\)
Since the order-$q$ moments are not all
$\pi$-a.s. constant, $H_{2q}\not\equiv0$. Consequently,
\[
\int_{B_1}H_{2q}(u)\,\mathrm du>0.
\]

Since the Taylor remainder is uniform in direction,
\[
\int_{B_\varepsilon}S(k)\,\mathrm dk
=
\varepsilon^{d+2q}
\left(
\int_{B_1}H_{2q}(u)\,\mathrm du+o(1)
\right).
\]
For $\varepsilon<2\pi$, there are no Bragg peaks in
$B_\varepsilon$. Hence
\[
\frac{\widehat\beta_{\widetilde\Phi}(B_\varepsilon)}
{\varepsilon^d}
=
\varepsilon^{2q}
\left(
\int_{B_1}H_{2q}(u)\,\mathrm du+o(1)
\right).
\]
The coefficient is strictly positive because
$H_{2q}$ is nonnegative, continuous, and not identically zero.
Thus, the process is $2q$-uniform but not beyond
$2q$-uniform.
\qed
\end{proo}


\section*{Acknowledgement}
Both authors are grateful to Prof. Dr. Günter Last for his valuable insights during the preparation of this article, which contributed significantly to shaping the results.

The work of Daniela Flimmel was supported by the Czech Science Foundation, project no. 26-22373N. The work of Luca Lotz was supported by the Deutsche Forschungsgemeinschaft (DFG, German Research Foundation) through the SPP 2265 under grant number KL 3391/2-2.

\section*{Statement on AI use}
The authors declare that they used AI tools for the following tasks: checking language and style, searching for relevant literature (the relevance of which the authors verified), and assisting in writing the code used to produce the figures. Specifically, the authors declare that they did not use AI tools to generate the scientific content.

\end{document}